\documentclass[10pt]{article}
\usepackage{amsmath,amssymb}
\usepackage{pstricks}
\usepackage{pst-node}
\usepackage[dvips]{graphicx}
\usepackage{pst-poly}
\usepackage{multido}

\newtheorem{result}{Theorem}

\renewcommand*{\thetheorem}{\Alph{theorem}}

\newcommand{\qed}{%
\ifmmode % if math mode, assume display: omit penalty etc.
\else \leavevmode\unskip\penalty9999 \hbox{}\nobreak\hfill \fi
\quad\hbox{\qedsymbol}}
\newcommand{\openbox}{\leavevmode \hbox to.77778em{%
\hfil\vrule
\vbox to.675em{\hrule width.6em\vfil\hrule}%
\vrule\hfil}}
\newcommand{\qedsymbol}{\openbox}

\newcommand{\showgrid}{}
\newcommand{\gridon}{\renewcommand{\showgrid}{\psset{subgriddiv=1,griddots=10,gridlabels=6pt}\psgrid}}
\gridon

\begin{document}
\begin{center}
{\bf\LARGE  A Proof Technique for Skewness of Graphs}
%Skewness of  Generalized Petersen Graphs and Related Graphs, II  }
\end{center}

\vskip8pt

\centerline{\large  Gek L. Chia$^{a,c}$  \   Chan L. Lee$^b$ \ and \ Yan Hao Ling$^b$ }

\begin{center}
\itshape\small $^a\/$Department of Mathematical and Actuarial Sciences, \\ Lee Kong Chian Faculty of Engineering and Science,  \\ Universiti Tunku Abdul Rahman, Sungai Long Campus,  \\  Cheras 43000 Kajang, Selangor,  Malaysia \\
\vspace{1mm}
$^b\/$Department of Mathematics, Statistics and Computing, \\ NUS High School of Mathematics \& Science,\\  20 Clementi Avenue 1, Singapore, 129957  \\ \vspace{1mm}
 $^c\/$Institute of Mathematical Sciences, University of Malaya, \\ 50603 Kuala Lumpur,  Malaysia  \\
\end{center}

\begin{abstract}
The {\em skewness\/} of a graph $G\/$ is the minimum number of edges in $G\/$ whose removal results in a planar graph.
By appropriately introducing a weight  to each edge of a graph, we  determine, among other thing, the skewness of the generalized Petersen graph $P(4k, k)\/$ for odd $k \geq 9\/$.
This provides an answer to the conjecture raised in \cite{cl3:refer}.
\end{abstract}

\vspace{3mm}

%\section{Introduction}

Let $G\/$ be a graph.  The {\em skewness\/} of $G\/$, denoted $\mu(G)\/$, is defined to be the minimum number of edges in $G\/$ whose removal results in a planar graph. Skewness of graph was first  introduced in the 1970's (see  \cite{kai:refer} - \cite{kai2:refer} and   \cite{pah:refer}). It was further explored in \cite{cim:refer}, \cite{lie:refer} and \cite{lg:refer}. More about the skewness of a graph can be found in \cite{cl2:refer} and \cite{cs:refer}.

%The skewness of a graph is clearly a lower bound of its crossing number. More about the skewness of a graph can be found in \cite{cl2:refer} and \cite{cs:refer}.

\vspace{1mm}
Let $n\/$ and $k\/$ be two integers such that $1 \leq k \leq n-1\/$. Recall that the {\em generalized Petersen graph\/} $P(n,k)\/$ is defined to have  vertex-set $\{u_i, v_i : i = 0, 1, \ldots, n-1\}\/$ and edge-set $\{ u_iu_{i+1}, u_iv_i, v_iv_{i+k} : i = 0, 1, \ldots, n-1\/$ with subscripts reduced modulo $n \}.\/$ % Edges of the form $u_iv_i\/$ are called the {\em spokes\/} of $P(n,k)\/$.

\vspace{1mm}
 Earlier, the authors in \cite{cl3:refer} showed that $\mu(P(4k, k)) = k + 1\/$ if $k \geq 4\/$ is even.  In the same paper, they conjectured that $\mu(P(4k, k)) = k + 2\/$  if $k \geq 5\/$ is odd.   We shall prove that the conjecture is true for $k \geq 9\/$.

\vspace{1mm} We shall first describe a family of graphs which is very much related to the generalized Petersen graph $P(sk, k)\/$.

\vspace{1mm}
Let $k\/$ and $s\/$ be two integers such that $k \geq 1\/$ and $s \geq 3\/$. The graph $Q_s(k)\/$ is defined to have vertex-set $\{ 0, 1, \ldots, sk-1, x_0, x_1, \ldots, x_{k-1}\}\/$ and edge-set $\{ i(i+1), jx_j, (j+k)x_j, (j+2k)x_j, \dots, (j+(s-1)k)x_j : i = 0, 1, \ldots, sk-1\/$, $j = 0, 1, \ldots, k-1\/$ with the operations reduced modulo $\/sk\/$ and those on the  subscripts reduced modulo $k \}.\/$ \

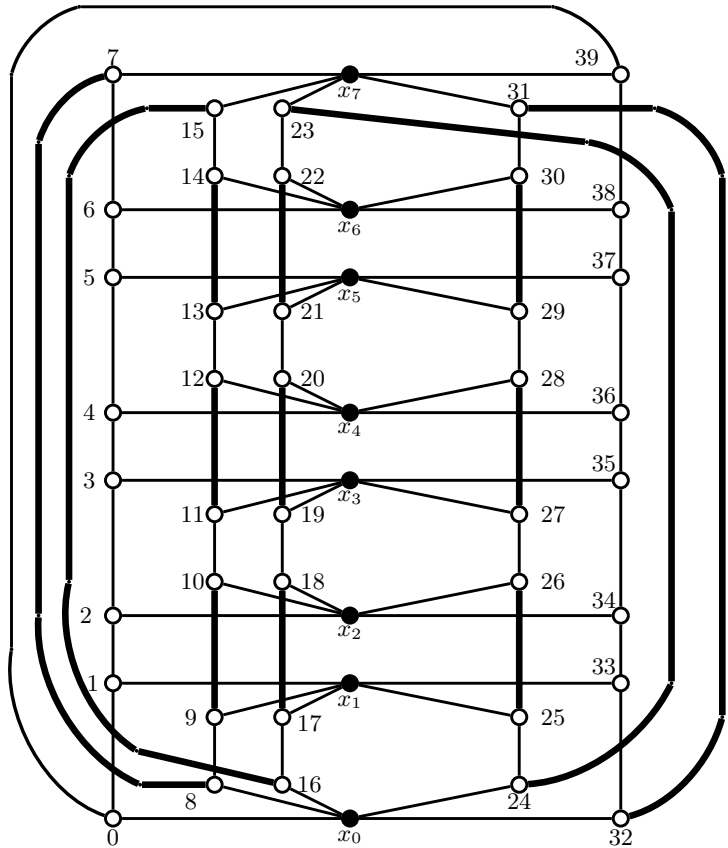
\begin{figure}[h]
\begin{center}
\psset{unit=0.45cm}
\begin{pspicture}(28,19)(1.5,-6)
   %\showgrid

\cnode*[linewidth=1.1pt](14,0){0.12cm}{x0}
\cnode*[linewidth=1.1pt](14,4){0.12cm}{x1}
\cnode*[linewidth=1.1pt](14,6){0.12cm}{x-4}
\cnode*[linewidth=1.1pt](14,10){0.12cm}{x-3}
\cnode*[linewidth=1.1pt](14,12){0.12cm}{x-2}
\cnode*[linewidth=1.1pt](14,16){0.12cm}{x-1}

\cnode[linewidth=1.1pt](7,0){0.12cm}{a0}
\cnode[linewidth=1.1pt](7,4){0.12cm}{a1}
\cnode[linewidth=1.1pt](7,6){0.12cm}{a-4}
\cnode[linewidth=1.1pt](7,10){0.12cm}{a-3}
\cnode[linewidth=1.1pt](7,12){0.12cm}{a-2}
\cnode[linewidth=1.1pt](7,16){0.12cm}{a-1}

\cnode[linewidth=1.1pt](22,0){0.12cm}{as}
\cnode[linewidth=1.1pt](22,4){0.12cm}{as1}
\cnode[linewidth=1.1pt](22,6){0.12cm}{as-4}
\cnode[linewidth=1.1pt](22,10){0.12cm}{as-3}
\cnode[linewidth=1.1pt](22,12){0.12cm}{as-2}
\cnode[linewidth=1.1pt](22,16){0.12cm}{as-1}

\ncline[linewidth=1.1pt]{-}{a1}{a-4} \ncline[linewidth=1.1pt]{-}{as-4}{as1}  %\ncline[linewidth=1.1pt]{-}{as1}{e3} \ncline[linewidth=1.1pt]{-}{as-4}{e4}

\cnode*[linewidth=1.1pt](4,-1){0.025cm}{f1} \cnode*[linewidth=1.1pt](4,16){0.025cm}{f2}
\cnode*[linewidth=1.1pt](6,18){0.025cm}{f3} \cnode*[linewidth=1.1pt](20,18){0.025cm}{f4}

\ncline[linewidth=1.1pt]{-}{f1}{f2}  \ncline[linewidth=1.1pt]{-}{f3}{f4}
\ncarc[arcangle=30,linewidth=1.3pt,linecolor=black]{-}{f2}{f3} \ncarc[arcangle=30,linewidth=1.3pt,linecolor=black]{-}{f4}{as-1}

\ncline[linewidth=1.1pt]{-}{x0}{a0}  \ncline[linewidth=1.1pt]{-}{x0}{as}
\ncline[linewidth=1.1pt]{-}{x1}{a1}  \ncline[linewidth=1.1pt]{-}{x1}{as1}
\ncline[linewidth=1.1pt]{-}{x-4}{a-4}  \ncline[linewidth=1.1pt]{-}{x-4}{as-4}
\ncline[linewidth=1.1pt]{-}{x-3}{a-3}  \ncline[linewidth=1.1pt]{-}{x-3}{as-3}
\ncline[linewidth=1.1pt]{-}{x-2}{a-2}  \ncline[linewidth=1.1pt]{-}{x-2}{as-2}
\ncline[linewidth=1.1pt]{-}{x-1}{a-1}  \ncline[linewidth=1.1pt]{-}{x-1}{as-1}

\ncline[linewidth=1.1pt]{-}{a0}{a1}  \ncline[linewidth=1.1pt]{-}{a-4}{a-3}  \ncline[linewidth=1.1pt]{-}{a-3}{a-2}  \ncline[linewidth=1.1pt]{-}{a-2}{a-1}
\ncline[linewidth=1.1pt]{-}{as}{as1}  \ncline[linewidth=1.1pt]{-}{as-4}{as-3}  \ncline[linewidth=1.1pt]{-}{as-3}{as-2}  \ncline[linewidth=1.1pt]{-}{as-2}{as-1}

\cnode[linewidth=1.1pt](10,1){0.12cm}{ak} \cnode[linewidth=1.1pt](12,1){0.12cm}{a2k}   \cnode[linewidth=1.1pt](19,1){0.12cm}{as-2k}
\cnode[linewidth=1.1pt](10,3){0.12cm}{ak1} \cnode[linewidth=1.1pt](12,3){0.12cm}{a2k1}   \cnode[linewidth=1.1pt](19,3){0.12cm}{as-2k1}

\ncline[linewidth=1.1pt]{-}{x0}{ak}   \ncline[linewidth=1.1pt]{-}{x0}{a2k} \ncline[linewidth=1.1pt]{-}{x0}{as-2k}
\ncline[linewidth=1.1pt]{-}{x1}{ak1}   \ncline[linewidth=1.1pt]{-}{x1}{a2k1} \ncline[linewidth=1.1pt]{-}{x1}{as-2k1}
\ncline[linewidth=1.1pt]{-}{ak1}{ak}   \ncline[linewidth=1.1pt]{-}{a2k1}{a2k} \ncline[linewidth=1.1pt]{-}{as-2k1}{as-2k}

\cnode[linewidth=1.1pt](10,7){0.12cm}{a2k-4} \cnode[linewidth=1.1pt](12,7){0.12cm}{a3k-4}   \cnode[linewidth=1.1pt](19,7){0.12cm}{ask-4}
\cnode[linewidth=1.1pt](10,9){0.12cm}{a2k-3} \cnode[linewidth=1.1pt](12,9){0.12cm}{a3k-3}   \cnode[linewidth=1.1pt](19,9){0.12cm}{ask-3}

\ncline[linewidth=1.1pt]{-}{x-4}{a2k-4}   \ncline[linewidth=1.1pt]{-}{x-4}{a3k-4} \ncline[linewidth=1.1pt]{-}{x-4}{ask-4}
\ncline[linewidth=1.1pt]{-}{x-3}{a2k-3}   \ncline[linewidth=1.1pt]{-}{x-3}{a3k-3} \ncline[linewidth=1.1pt]{-}{x-3}{ask-3}
\ncline[linewidth=1.1pt]{-}{a2k-3}{a2k-4}   \ncline[linewidth=1.1pt]{-}{a3k-3}{a3k-4} \ncline[linewidth=1.1pt]{-}{ask-3}{ask-4}

\cnode[linewidth=1.1pt](10,13){0.12cm}{a2k-2} \cnode[linewidth=1.1pt](12,13){0.12cm}{a3k-2}   \cnode[linewidth=1.1pt](19,13){0.12cm}{ask-2}
\cnode[linewidth=1.1pt](10,15){0.12cm}{a2k-1} \cnode[linewidth=1.1pt](12,15){0.12cm}{a3k-1}   \cnode[linewidth=1.1pt](19,15){0.12cm}{ask-1}

\ncline[linewidth=1.1pt]{-}{x-2}{a2k-2}   \ncline[linewidth=1.1pt]{-}{x-2}{a3k-2} \ncline[linewidth=1.1pt]{-}{x-2}{ask-2}
\ncline[linewidth=1.1pt]{-}{x-1}{a2k-1}   \ncline[linewidth=1.1pt]{-}{x-1}{a3k-1} \ncline[linewidth=1.1pt]{-}{x-1}{ask-1}
\ncline[linewidth=1.1pt]{-}{a2k-2}{a2k-1}   \ncline[linewidth=1.1pt]{-}{a3k-2}{a3k-1} \ncline[linewidth=1.1pt]{-}{ask-2}{ask-1}

%\rput(14,-8){ $Q_5(8)\/$}

\rput(6.4,-2){\small $1$}  \rput(6.2,0){\small $2$}   \rput(14,-.5){\small$x_2$}   \rput(14,-2.55){\small$x_1$}    %\rput(22,1.45){\small $(s-1)k$}
\rput(6.3,4){\small $3$}   \rput(14,3.5){\small$x_3$}  \rput(20,3){\small $27$}    \rput(21.5,-1.5){\small $33$}     \rput(21.5,4.5){\small $35$}
\rput(6.3,6){\small $4$}   \rput(14,5.5){\small$x_{4}$}  \rput(21.5,6.5){\small $36$}
\rput(6.3,10){\small $5$}   \rput(14,9.45){\small$x_{5}$}  \rput(21.5,10.5){\small $37$}
\rput(6.3,12){\small $6$}   \rput(14,11.45){\small$x_6$}  \rput(21.5,12.5){\small $38$}
\rput(7,16.5){\small $7$}   \rput(14,15.45){\small$x_7$}  \rput(21,16.5){\small $39$}

\rput(9.35,1){\small $10$} \rput(12.9,1){\small $18$}
\rput(9.35,3){\small $11$} \rput(12.9,3){\small $19$}  \rput(20,1){\small $26$}

\rput(9.35,7){\small $12$} \rput(12.9,7){\small $20$}  \rput(20,7){\small $28$}
\rput(9.35,9){\small $13$} \rput(12.9,9){\small $21$}  \rput(20,9){\small $29$}

\rput(9.35,13){\small $14$} \rput(12.9,13){\small $22$}  \rput(19,15.5){\small $31$}
\rput(9.35,14.35){\small $15$} \rput(12.6,14.35){\small $23$}  \rput(20,13){\small $30$}

%extension

\cnode*[linewidth=1.1pt](14,-6){0.12cm}{-x0}  \cnode[linewidth=1.1pt](7,-6){0.12cm}{-a0}   \cnode[linewidth=1.1pt](22,-6){0.12cm}{-as}

\cnode*[linewidth=1.1pt](14,-2){0.12cm}{x00}

\cnode[linewidth=1.1pt](10,-3){0.12cm}{-ak} \cnode[linewidth=1.1pt](12,-3){0.12cm}{-a2k}   \cnode[linewidth=1.1pt](19,-3){0.12cm}{-as-2k}
\cnode[linewidth=1.1pt](10,-5){0.12cm}{-ak1} \cnode[linewidth=1.1pt](12,-5){0.12cm}{-a2k1}   \cnode[linewidth=1.1pt](19,-5){0.12cm}{-as-2k1}

\ncline[linewidth=1.1pt]{-}{x00}{-ak}   \ncline[linewidth=1.1pt]{-}{x00}{-a2k} \ncline[linewidth=1.1pt]{-}{x00}{-as-2k}
\ncline[linewidth=1.1pt]{-}{-x0}{-ak1}   \ncline[linewidth=1.1pt]{-}{-x0}{-a2k1} \ncline[linewidth=1.1pt]{-}{-x0}{-as-2k1}
\ncline[linewidth=1.1pt]{-}{-ak1}{-ak}   \ncline[linewidth=1.1pt]{-}{-a2k1}{-a2k} \ncline[linewidth=1.1pt]{-}{-as-2k1}{-as-2k}

\ncline[linewidth=1.1pt]{-}{x00}{a00}   \ncline[linewidth=1.1pt]{-}{x00}{ass}

\ncline[linewidth=1.1pt]{-}{a0}{-a0}  \ncline[linewidth=1.1pt]{-}{as}{-as}

\ncline[linewidth=1.1pt]{-}{-x0}{-a0}  \ncline[linewidth=1.1pt]{-}{-x0}{-as}

 \cnode[fillstyle=solid,fillcolor=white,linewidth=1.1pt](7,-2){0.12cm}{a00}
 \cnode[fillstyle=solid,fillcolor=white,linewidth=1.1pt](22,-2){0.12cm}{ass}

\ncline[linewidth=1.1pt]{-}{x00}{a00}   \ncline[linewidth=1.1pt]{-}{x00}{ass}

\rput(7,-6.55){\small $0$}   \rput(14,-6.55){\small$x_0$}  \rput(22,-6.55){\small $32$}    \rput(21.5,.45){\small $34$}

\rput(9.3,-5.5){\small $8$} \rput(12.8,-5){\small $16$}  \rput(19,-5.5){\small $24$}
\rput(9.3,-3){\small $9$} \rput(12.8,-3.2){\small $17$}  \rput(20,-3){\small $25$}

\ncarc[arcangle=-40,linewidth=1.3pt,linecolor=black]{-}{f1}{-a0}

%thick lines

\ncline[linewidth=2.5pt]{-}{-ak}{ak}  \ncline[linewidth=2.5pt]{-}{a2k-4}{ak1}     \ncline[linewidth=2.5pt]{-}{a2k-3}{a2k-2}
\ncline[linewidth=2.5pt]{-}{-a2k}{a2k}  \ncline[linewidth=2.5pt]{-}{a3k-4}{a2k1}     \ncline[linewidth=2.5pt]{-}{a3k-3}{a3k-2}
\ncline[linewidth=2.5pt]{-}{-as-2k}{as-2k}  \ncline[linewidth=2.5pt]{-}{ask-4}{as-2k1}     \ncline[linewidth=2.5pt]{-}{ask-3}{ask-2}

\cnode*[linewidth=1.1pt](4.8,0){0.025cm}{g1} \cnode*[linewidth=1.1pt](4.8,14){0.025cm}{g2}   \cnode*[linewidth=1.5pt](7.8,-5){0.025cm}{g3}
\ncline[linewidth=2.5pt]{-}{g1}{g2}  \ncline[linewidth=2.5pt]{-}{-ak1}{g3}
\ncarc[arcangle=-35,linewidth=2.5pt,linecolor=black]{-}{g1}{g3}     \ncarc[arcangle=30,linewidth=2.5pt,linecolor=black]{-}{g2}{a-1}

\cnode*[linewidth=1.1pt](5.7,1){0.025cm}{h1} \cnode*[linewidth=1.1pt](5.7,13){0.025cm}{h2}   \cnode*[linewidth=1.5pt](7.7,-4){0.025cm}{h3}
\cnode*[linewidth=1.5pt](8,15){0.025cm}{h4}

\ncline[linewidth=2.5pt]{-}{h1}{h2}    \ncline[linewidth=2.5pt]{-}{-a2k1}{h3}    \ncline[linewidth=2.5pt]{-}{a2k-1}{h4}
\ncarc[arcangle=-35,linewidth=2.5pt,linecolor=black]{-}{h1}{h3}    \ncarc[arcangle=30,linewidth=2.5pt,linecolor=black]{-}{h2}{h4}

\cnode*[linewidth=1.1pt](23,15){0.025cm}{k1} \cnode*[linewidth=1.1pt](25,13){0.025cm}{k2}   \cnode*[linewidth=1.5pt](25,-3){0.025cm}{k3}
\ncline[linewidth=2.5pt]{-}{k3}{k2}    \ncline[linewidth=2.5pt]{-}{k1}{ask-1}
 \ncarc[arcangle=30,linewidth=2.5pt,linecolor=black]{-}{k1}{k2}  \ncarc[arcangle=30,linewidth=2.5pt,linecolor=black]{-}{k3}{-as}

\cnode*[linewidth=1.1pt](21,14){0.025cm}{l1} \cnode*[linewidth=1.1pt](23.5,12){0.025cm}{l2}   \cnode*[linewidth=1.5pt](23.5,-2){0.025cm}{l3}
\ncline[linewidth=2.5pt]{-}{l3}{l2}   \ncline[linewidth=2.5pt]{-}{l1}{a3k-1}
\ncarc[arcangle=30,linewidth=2.5pt,linecolor=black]{-}{l1}{l2}   \ncarc[arcangle=30,linewidth=2.5pt,linecolor=black]{-}{l3}{-as-2k1}

\end{pspicture}
\end{center}
\caption{A drawing of the graph $Q_5(8)\/$.}  \label{Q5(8)}
\end{figure}

\vspace{2mm}
Drawings for the graphs $Q_3(k)\/$ and $Q_4(k)\/$ can be found in the papers \cite{cl:refer} and \cite{cl3:refer} respectively. A drawing for the graph $Q_5(8)\/$ is depicted in Figure \ref{Q5(8)}.

\vspace{2mm} Throughout this paper, if $\/G\/$ has skewness $\/r\/$, then we let ${\mathcal R} (G)\/$ denote a set of $\/r\/$ edges in $\/G\/$ whose removal results in a planar graph.

\vspace{2mm}
In \cite{cl:refer} and \cite{cl3:refer} it was shown that $\mu(Q_3(k))= \lceil k/2\rceil +1 $ and $\mu(Q_4(k))= k +1 $ respectively. A more general result is now established by employing a proof technique that has not been used before.

\vspace{2mm}
Suppose $G\/$ is a weighted graph and $e\/$ is an edge of $G\/$ with weight $w(e)\/$.  If $H\/$ is any subgraph (proper or improper) of $G\/$, we let $W(H) = \sum_{e \in E(H)} w(e)\/$ denote the {\em weight of $H\/$\/}. In particular, if $G\/$ is a plane graph and $F\/$ is an face of $G\/$, we let $W(F) = \sum_{e \in E(F)} w(e)\/$ denote the {\em weight of $F\/$\/} in $G\/$.

\vspace{2mm}
\begin{result}  \label{skq4}
  $\mu(Q_s(k))= \lceil (s-2)k/2\rceil +1 $  \ if \  $k\ge 4.$
\end{result}

\vspace{1mm} \noindent
{\sc Proof:} Let $e\/$ be an edge of $Q_s(k)\/$ and let
\[       w(e)  = \left\{ \begin{array}{ll}
              2  & \ \ \ \ \  \ \mbox{if \ $e = i(i+1)$,}   \\
                 &      \\
            k-2  & \ \ \ \ \ \ \mbox{otherwise }
                                     \end{array}
                 \right.   \]
for $i=0, 1, \ldots, sk-1\/$.  Then $W(Q_s(k)) = sk^2\/$.

\vspace{2mm}
Hence, if $J\/$ is a graph obtained from $Q_s(k)\/$ by deleting a  set of $t\/$ edges, then $W(J) \leq sk^2 -2t\/$ (because $k-2 \geq 2\/$).

 \vspace{2mm}
 Let $H\/$ denote a planar graph obtained from $Q_s(k)\/$ by deleting a set of $t'= \mu(Q_s(k))\/$ edges. Then the number of faces in $H\/$ is $sk-k-t'+2\/$. It is easy to see that $W(F) \geq 4k-4\/$ for any face $F\/$ in $H\/$.
  Since each edge is contained  in $2\/$ faces, we have the following inequality
 \[   2(sk^2 -2t') \geq (4k-4)(sk-k-t'+2)    \]
which gives $t' \geq (s-2)k/2 +1\/$. This proves the lower bound.

\vspace{2mm}
To prove the upper bound, we show the existence of a spanning planar subgraph $H_s(k)\/$ obtained by deleting a set of $\lceil (s-2)k/2\rceil +1 $  edges   from $Q_s(k)\/$.

\vspace{2mm}
Figures  2, 3 and 4  depicts three drawings of $H_s(k)\/$ according to the parities on $s\/$ and $k\/$. When $k\/$ is even, the set of  edges that have been deleted from $Q_s(k)\/$ is  given by $(2i-1)(2i)\/$, $i= k/2, k/2+1, \ldots, k/2 + k(s-2)/2\/$ (see Figure 2).   Note also that if the $9\/$ thick edges in the graph $Q_5(8)\/$ (as depicted in Figure \ref{Q5(8)}) are deleted, we obtain the graph $H_5(8)\/$.

\vspace{2mm}
When $k\/$ is odd, the set of  edges that have been deleted from $Q_s(k)\/$ is  given by $(sk-1)0$, $(2i-1)(2i)\/$, $i= (k-1)/2+1, (k-1)/2+2, \ldots, (k-1)/2 + \lceil k(s-2)/2\rceil\/$ (see Figures 3 and 4).

\vspace{2mm}  This completes the proof.   \qed

\vspace{3mm}

\begin{result}
$\mu(P(4k,k)) = k+2\/$ if $k \geq 9\/$ is odd.
\end{result}

\vspace{1mm} \noindent
{\sc Proof:} The upper bound $\mu(P(4k,k)) \leq  k+2\/$ for odd  $k \geq 5\/$ was established in \cite{cl3:refer}. Hence we just need to show that $\mu(P(4k,k)) \geq  k+2\/$ for odd  $k \geq 9\/$.

\vspace{2mm}
First we shall show that $\mu(P(4k,k)) \geq  k+1\/$ if  $k \geq 9\/$ is odd.  The proof presented here employed a new technique (which involves assigning appropriate weights to its edges) and hence  is different from the one given in \cite{cl3:refer}.

\vspace{2mm}
Let $e\/$ be an edge of $P(4k,k)\/$ and let
\[       w(e)  = \left\{ \begin{array}{ll}
              4  & \ \ \ \ \  \ \mbox{if \ $e = u_iu_{i+1}$,}   \\
                 &      \\
            k-3  & \ \ \ \ \ \ \mbox{if \ $e= u_iv_i,$} \\
                 &      \\
            2k-2 & \ \ \ \ \ \ \mbox{if \ $e=v_iv_{i+k}$ }
                                     \end{array}
                 \right.   \]
for $i=0, 1, \ldots, n-1\/$. Then $W(P(4k,k)) = 4k(3k-1)\/$.        Moreover,   it is easy to see that,  for any cycle $C\/$ in $P(4k,k)\/$,  $W(C) \geq 8k-8\/$ and that equality holds if and only if $C\/$ is any of the following types.

\vspace{1mm}
(i) $u_iu_{i+1}u_{i+2} \ldots u_{i+k-1}u_{i+k}v_{i+k}v_iu_i\/$,

\vspace{1mm}
(ii) $u_iu_{i+1}v_{i+1}v_{i+k+1}u_{i+k+1}u_{i+k}v_{i+k}v_iu_i\/$,

\vspace{1mm} (iii) $u_iu_{i+1}v_{i+1}v_{i-k+1}u_{i-k+1}u_{i-k}v_{i-k}v_iu_i\/$,

\vspace{1mm}
(iv) $v_iv_{i+k}v_{i+2k}v_{i+3k}v_i\/$.

\vspace{2mm}
Let $H\/$ denote a planar graph obtained by deleting a set of $t=\mu(P(4k,k))\/$ edges from $P(4k,k)\/$. Then $W(H) \leq 4k(3k-1) - 4t\/$ (since $k \geq 9\/$ implies that  $k-3 > 4\/$ and $2k-2 > 4$).

\vspace{2mm}
From Euler's formula for plane graph, we see that  the number of faces in $H\/$ is $4k-t+2\/$. Since each edge is contained in only 2 faces, we have the following inequality
\[    2(4k(3k-1) -4t) \geq (4k-t+2)(8k-8)  \]
which gives $t \geq k+1\/$.

\vspace{2mm}
Now suppose $\mu(P(4k,k)) = k+1\/$. The fact that equality is tight implies the following.

\vspace{2mm} (a) Only faces of the types (i) to (iv) are found in the planar graph $H\/$.  (b) The tight equality $4k(3k-1) -4t\/$ implies that ${\mathcal R} (P(4k,k))\/$ consists of edges of the form $u_iu_{i+1}\/$. (c) Since $P(4k,k)\/$ is a $3\/$-regular graph, removing any two adjacent edges yields a pendant vertex, resulting in  $W(F) > 8k-8\/$. Hence we conclude that ${\mathcal R} (P(4k,k))\/$ consists of only independent edges of the form $u_iu_{i+1}\/$.

\vspace{2mm}
Now for any edge $e\/$ of $P(4k,k)\/$, let
\[       w'(e)  = \left\{ \begin{array}{ll}
              0  & \ \ \ \ \  \ \mbox{if \ $e = u_iu_{i+1}$,}   \\
                 &      \\
              1  & \ \ \ \ \ \ \mbox{if \ $e= u_iv_i,$} \\
                 &      \\
              2  & \ \ \ \ \ \ \mbox{if \ $e=v_iv_{i+k}$ }
                                     \end{array}
                 \right.   \]
for $i=0, 1, \ldots, n-1\/$. We note that, if $F\/$ is a face of $H\/$, then $W'(F) = \sum_{e \in E(F)} w'(e)\/$ is equal to $4\/$ if $F\/$ is of type (i) and $W'(F)\/$ is equal to $8\/$ if $F\/$ is of type (ii), (iii) or (iv).

\vspace{2mm}
Let $x\/$ be the number of faces of type (i) in $H\/$ and let $\cal F _H\/$ denote the set of all faces in $H\/$. Then $\sum_{F \in \cal F _H} W'(F) = 4x +8(3k+1-x)\/$.  Since ${\mathcal R} (P(4k,k))\/$ consists of only independent edges of the form $u_iu_{i+1}\/$, we see that $W'(H) = 12k\/$. As such, we have
\[    4x +8(3k+1-x) = 2 (12k)  \]
which gives $x=2\/$.

\vspace{2mm}
Now let the two faces of type (i) be

\vspace{1mm}
$F_1 = u_yu_{y+1} u_{y+2} \ldots u_{y+k-1}u_{y+k}v_{y+k}v_yu_y\/$  \  \ and

\vspace{1mm} $F_2 = u_zu_{z+1} u_{z+2} \ldots u_{z+k-1}u_{z+k}v_{z+k}v_zu_z\/$.

\vspace{2mm}
Let $A = \{y+1, y+2, \ldots, y+k-1\} \cap \{ z+1, z+2, \ldots, y+k-1\}\/$.  We assert that $A = \emptyset\/$.

\vspace{2mm} Suppose $h \in A\/$. Then the edges $u_{h-1}u_h, u_hu_{h+1}\/$ are in $F_1  \cap F_2 \/$ implying that $u_hv_h \in  {\mathcal R} (P(4k,k))\/$, a contradiction.

\vspace{2mm}
Call a vertex $u_i\/$ of $P(4k,k)\/$ a {\em good\/} vertex if all edges incident to it are not in  ${\mathcal R} (P(4k,k))\/$; otherwise it is called a {\em bad\/} vertex. Since ${\mathcal R} (P(4k,k))\/$ consists of independent edges of the form $u_iu_{i+1}\/$, the number of bad vertices is $2\mu(P(4k,k))\/$.

\vspace{2mm}
Since $y+1, y+2, \ldots, y+k-1,  z+1, z+2, \ldots, y+k-1\/$ are distinct, and since $F_1\/$ and $F_2\/$ share at most one common edge (of the form $u_uv_i\/$) in the planar graph $H\/$,  the number of good vertices in $H\/$ is at least $2k-2\/$ (given by the distinct vertices $u_{y+1}, u_{y+2}, \ldots, u_{y+k-1},  u_{z+1},   u_{z+2},  \ldots, u_{z+k-1}\/$) and at most $2k-1\/$ (if $F_1\/$ and $F_2\/$ have an edge in common). Hence $H\/$ has at least $2k+1\/$ bad vertices which are either consecutive  vertices of the form $u_n, u_{n+1}, \ldots,  u_{n+2k-1}\/$ or are separated into $2\/$ sets of consecutive vertices by $V(F_1) \cup V(F_2)\/$. In any case, we can find at least $k+1\/$  consecutive bad vertices. Since  we can relabel the vertices if necessary, we may assume without loss of generality that these vertices are $u_0, u_1, \ldots, u_k\/$, and that ${\mathcal R} (P(4k,k))\/$ contains $u_0u_1, u_2u_3, \ldots, u_{k-1}u_k\/$.

\vspace{2mm}
Let $J\/$ denote the subgraph obtained from $P(4k,k)\/$ by deleting $u_1, u_2, \ldots,  \linebreak u_{k-1}, v_1, v_2, \ldots, v_{k-1}\/$ together with the edge $v_{3k}v_k\/$. Then $J\/$ is a subdivision of $Q_3(k)\/$. Note that the vertices of degree-$2\/$ in $J\/$ are $u_0, v_0, u_k, v_{k+1}, v_{k+2}, \linebreak \ldots,  v_{2k-1}, v_{4k-1}, v_{4k-2}, \ldots, v_{3k}\/$. When these degree-$2\/$ vertices are suppress\-ed, the resulting graph is isomorphic to $Q_3(k)\/$ with $v_ku_{k+1}u_{k+2} \ldots u_{4k-1}v_k\/$ playing role of the $3k\/$-cycle in $Q_3(k)\/$.

\vspace{2mm}
Recall that $\mu(Q_3(k)) = (k+3)/2\/$. The preceding arguments imply that
\[ \mu(P(4k,k)) \geq |\{ u_0u_1, u_2u_3, \ldots, u_{k-1}u_k \}| +    \mu(Q_3(k))  = k+2    \]
which is a contradiction.   Hence $\mu(P(4k,k)) = k+2\/$.                \hfill $\Box$

\begin{figure}[ht]
\begin{center}
\psset{unit=0.45cm}
\begin{pspicture}(28,30)(2,-6)
   %\showgrid

\cnode*[linewidth=1.1pt](14,2){0.12cm}{x0}
\cnode*[linewidth=1.1pt](14,8){0.12cm}{x1}
\cnode*[linewidth=1.1pt](14,14){0.12cm}{x-4}
\cnode*[linewidth=1.1pt](14,20){0.12cm}{x-3}
\cnode*[linewidth=1.1pt](14,22){0.12cm}{x-2}
\cnode*[linewidth=1.1pt](14,28){0.12cm}{x-1}

\cnode[linewidth=1.1pt](7,2){0.12cm}{a0}
\cnode[linewidth=1.1pt](7,8){0.12cm}{a1}
\cnode[linewidth=1.1pt](7,14){0.12cm}{a-4}
\cnode[linewidth=1.1pt](7,20){0.12cm}{a-3}
\cnode[linewidth=1.1pt](7,22){0.12cm}{a-2}
\cnode[linewidth=1.1pt](7,28){0.12cm}{a-1}

\cnode[linewidth=1.1pt](22,2){0.12cm}{as}
\cnode[linewidth=1.1pt](22,8){0.12cm}{as1}
\cnode[linewidth=1.1pt](22,14){0.12cm}{as-4}
\cnode[linewidth=1.1pt](22,20){0.12cm}{as-3}
\cnode[linewidth=1.1pt](22,22){0.12cm}{as-2}
\cnode[linewidth=1.1pt](22,28){0.12cm}{as-1}

\cnode*[linewidth=1.1pt](7,10){0.05cm}{d0} \cnode*[linewidth=1.1pt](7,11){0.05cm}{d1} \cnode*[linewidth=1.1pt](7,12){0.05cm}{d2}
\cnode*[linewidth=1.1pt](22,10){0.05cm}{d3} \cnode*[linewidth=1.1pt](22,11){0.05cm}{d4} \cnode*[linewidth=1.1pt](22,12){0.05cm}{d5}
\cnode*[linewidth=1.1pt](14,10){0.05cm}{d6} \cnode*[linewidth=1.1pt](14,11){0.05cm}{d7} \cnode*[linewidth=1.1pt](14,12){0.05cm}{d8}

\cnode*[linewidth=1.1pt](14,17){0.05cm}{d9} \cnode*[linewidth=1.1pt](15.5,17){0.05cm}{d10} \cnode*[linewidth=1.1pt](17,17){0.05cm}{d11}
\cnode*[linewidth=1.1pt](14,5){0.05cm}{d9} \cnode*[linewidth=1.1pt](15.5,5){0.05cm}{d10} \cnode*[linewidth=1.1pt](17,5){0.05cm}{d11}
\cnode*[linewidth=1.1pt](14,25){0.05cm}{d9} \cnode*[linewidth=1.1pt](15.5,25){0.05cm}{d10} \cnode*[linewidth=1.1pt](17,25){0.05cm}{d11}

\cnode*[linewidth=1.1pt](7,9.2){0.025cm}{e1} \cnode*[linewidth=1.1pt](7,12.8){0.025cm}{e2}
\cnode*[linewidth=1.1pt](22,9.2){0.025cm}{e3} \cnode*[linewidth=1.1pt](22,12.8){0.025cm}{e4}
\ncline[linewidth=1.1pt]{-}{a1}{e1} \ncline[linewidth=1.1pt]{-}{a-4}{e2}  \ncline[linewidth=1.1pt]{-}{as1}{e3} \ncline[linewidth=1.1pt]{-}{as-4}{e4}

\cnode*[linewidth=1.1pt](4,-1){0.025cm}{f1} \cnode*[linewidth=1.1pt](4,28){0.025cm}{f2}
\cnode*[linewidth=1.1pt](6,30){0.025cm}{f3} \cnode*[linewidth=1.1pt](20,30){0.025cm}{f4}

\ncline[linewidth=1.1pt]{-}{f1}{f2}  \ncline[linewidth=1.1pt]{-}{f3}{f4}
 \ncarc[arcangle=30,linewidth=1.3pt,linecolor=black]{-}{f2}{f3} \ncarc[arcangle=30,linewidth=1.3pt,linecolor=black]{-}{f4}{as-1}

\ncline[linewidth=1.1pt]{-}{x0}{a0}  \ncline[linewidth=1.1pt]{-}{x0}{as}
\ncline[linewidth=1.1pt]{-}{x1}{a1}  \ncline[linewidth=1.1pt]{-}{x1}{as1}
\ncline[linewidth=1.1pt]{-}{x-4}{a-4}  \ncline[linewidth=1.1pt]{-}{x-4}{as-4}
\ncline[linewidth=1.1pt]{-}{x-3}{a-3}  \ncline[linewidth=1.1pt]{-}{x-3}{as-3}
\ncline[linewidth=1.1pt]{-}{x-2}{a-2}  \ncline[linewidth=1.1pt]{-}{x-2}{as-2}
\ncline[linewidth=1.1pt]{-}{x-1}{a-1}  \ncline[linewidth=1.1pt]{-}{x-1}{as-1}

\ncline[linewidth=1.1pt]{-}{a0}{a1}  \ncline[linewidth=1.1pt]{-}{a-4}{a-3}  \ncline[linewidth=1.1pt]{-}{a-3}{a-2}  \ncline[linewidth=1.1pt]{-}{a-2}{a-1}
\ncline[linewidth=1.1pt]{-}{as}{as1}  \ncline[linewidth=1.1pt]{-}{as-4}{as-3}  \ncline[linewidth=1.1pt]{-}{as-3}{as-2}  \ncline[linewidth=1.1pt]{-}{as-2}{as-1}

\cnode[linewidth=1.1pt](10,4){0.12cm}{ak} \cnode[linewidth=1.1pt](12,4){0.12cm}{a2k}   \cnode[linewidth=1.1pt](19,4){0.12cm}{as-2k}
\cnode[linewidth=1.1pt](10,6){0.12cm}{ak1} \cnode[linewidth=1.1pt](12,6){0.12cm}{a2k1}   \cnode[linewidth=1.1pt](19,6){0.12cm}{as-2k1}

\ncline[linewidth=1.1pt]{-}{x0}{ak}   \ncline[linewidth=1.1pt]{-}{x0}{a2k} \ncline[linewidth=1.1pt]{-}{x0}{as-2k}
\ncline[linewidth=1.1pt]{-}{x1}{ak1}   \ncline[linewidth=1.1pt]{-}{x1}{a2k1} \ncline[linewidth=1.1pt]{-}{x1}{as-2k1}
\ncline[linewidth=1.1pt]{-}{ak1}{ak}   \ncline[linewidth=1.1pt]{-}{a2k1}{a2k} \ncline[linewidth=1.1pt]{-}{as-2k1}{as-2k}

\cnode[linewidth=1.1pt](10,16){0.12cm}{a2k-4} \cnode[linewidth=1.1pt](12,16){0.12cm}{a3k-4}   \cnode[linewidth=1.1pt](19,16){0.12cm}{ask-4}
\cnode[linewidth=1.1pt](10,18){0.12cm}{a2k-3} \cnode[linewidth=1.1pt](12,18){0.12cm}{a3k-3}   \cnode[linewidth=1.1pt](19,18){0.12cm}{ask-3}

\ncline[linewidth=1.1pt]{-}{x-4}{a2k-4}   \ncline[linewidth=1.1pt]{-}{x-4}{a3k-4} \ncline[linewidth=1.1pt]{-}{x-4}{ask-4}
\ncline[linewidth=1.1pt]{-}{x-3}{a2k-3}   \ncline[linewidth=1.1pt]{-}{x-3}{a3k-3} \ncline[linewidth=1.1pt]{-}{x-3}{ask-3}
\ncline[linewidth=1.1pt]{-}{a2k-3}{a2k-4}   \ncline[linewidth=1.1pt]{-}{a3k-3}{a3k-4} \ncline[linewidth=1.1pt]{-}{ask-3}{ask-4}

\cnode[linewidth=1.1pt](10,24){0.12cm}{a2k-2} \cnode[linewidth=1.1pt](12,24){0.12cm}{a3k-2}   \cnode[linewidth=1.1pt](19,24){0.12cm}{ask-2}
\cnode[linewidth=1.1pt](10,26){0.12cm}{a2k-1} \cnode[linewidth=1.1pt](12,26){0.12cm}{a3k-1}   \cnode[linewidth=1.1pt](19,26){0.12cm}{ask-1}

\ncline[linewidth=1.1pt]{-}{x-2}{a2k-2}   \ncline[linewidth=1.1pt]{-}{x-2}{a3k-2} \ncline[linewidth=1.1pt]{-}{x-2}{ask-2}
\ncline[linewidth=1.1pt]{-}{x-1}{a2k-1}   \ncline[linewidth=1.1pt]{-}{x-1}{a3k-1} \ncline[linewidth=1.1pt]{-}{x-1}{ask-1}
\ncline[linewidth=1.1pt]{-}{a2k-2}{a2k-1}   \ncline[linewidth=1.1pt]{-}{a3k-2}{a3k-1} \ncline[linewidth=1.1pt]{-}{ask-2}{ask-1}

%\rput(14,-8){(a) $k\/$ is even.}

\rput(6.2,0){\small $1$}  \rput(6.2,2){\small $2$}   \rput(14,1.45){\small$x_2$}   \rput(14,-0.55){\small$x_1$}    %\rput(22,1.45){\small $(s-1)k$}
\rput(6.2,8){\small $3$}   \rput(14,8.55){\small$x_3$}  \rput(24.3,8){\small $(s-1)k+3$}    \rput(24.3,0){\small $(s-1)k+1$}
\rput(5.8,14){\small $k-4$}   \rput(14,13.45){\small$x_{k-4}$}  \rput(23.3,14){\small $sk-4$}
\rput(5.8,20){\small $k-3$}   \rput(14,20.5){\small$x_{k-3}$}  \rput(23.3,20){\small $sk-3$}
\rput(5.8,22){\small $k-2$}   \rput(14,21.5){\small$x_{k-2}$}  \rput(23.3,22){\small $sk-2$}
\rput(5.8,28){\small $k-1$}   \rput(14,28.5){\small$x_{k-1}$}  \rput(23.3,28){\small $sk-1$}

\rput(8.75,4){\small $k+2$} \rput(13.4,4){\small $2k+2$}  \rput(19.8,3.4){\small $(s-2)k+2$}
\rput(8.75,6){\small $k+3$} \rput(13.4,6){\small $2k+3$}  \rput(19.8,6.6){\small $(s-2)k+3$}

\rput(8.5,16){\small $2k-4$} \rput(13.6,16){\small $3k-4$}  \rput(19.8,15.3){\small $(s-1)k-4$}
\rput(8.5,18){\small $2k-3$} \rput(13.6,18){\small $3k-3$}  \rput(19.8,18.6){\small $(s-1)k-3$}

\rput(8.5,24){\small $2k-2$} \rput(13.6,24){\small $3k-2$}  \rput(19.8,23.3){\small $(s-1)k-2$}
\rput(8.5,26){\small $2k-1$} \rput(13.6,26){\small $3k-1$}  \rput(19.8,26.6){\small $(s-1)k-1$}

%extension

\cnode*[linewidth=1.1pt](14,-6){0.12cm}{-x0}  \cnode[linewidth=1.1pt](7,-6){0.12cm}{-a0}   \cnode[linewidth=1.1pt](22,-6){0.12cm}{-as}

\cnode*[linewidth=1.1pt](14,0){0.12cm}{x00}

\cnode*[linewidth=1.1pt](14,-3){0.05cm}{d9} \cnode*[linewidth=1.1pt](15.5,-3){0.05cm}{d10} \cnode*[linewidth=1.1pt](17,-3){0.05cm}{d11}

\cnode[linewidth=1.1pt](10,-2){0.12cm}{-ak} \cnode[linewidth=1.1pt](12,-2){0.12cm}{-a2k}   \cnode[linewidth=1.1pt](19,-2){0.12cm}{-as-2k}
\cnode[linewidth=1.1pt](10,-4){0.12cm}{-ak1} \cnode[linewidth=1.1pt](12,-4){0.12cm}{-a2k1}   \cnode[linewidth=1.1pt](19,-4){0.12cm}{-as-2k1}

\ncline[linewidth=1.1pt]{-}{x00}{-ak}   \ncline[linewidth=1.1pt]{-}{x00}{-a2k} \ncline[linewidth=1.1pt]{-}{x00}{-as-2k}
\ncline[linewidth=1.1pt]{-}{-x0}{-ak1}   \ncline[linewidth=1.1pt]{-}{-x0}{-a2k1} \ncline[linewidth=1.1pt]{-}{-x0}{-as-2k1}
\ncline[linewidth=1.1pt]{-}{-ak1}{-ak}   \ncline[linewidth=1.1pt]{-}{-a2k1}{-a2k} \ncline[linewidth=1.1pt]{-}{-as-2k1}{-as-2k}

\ncline[linewidth=1.1pt]{-}{x00}{a00}   \ncline[linewidth=1.1pt]{-}{x00}{ass}

\ncline[linewidth=1.1pt]{-}{a0}{-a0}  \ncline[linewidth=1.1pt]{-}{as}{-as}

\ncline[linewidth=1.1pt]{-}{-x0}{-a0}  \ncline[linewidth=1.1pt]{-}{-x0}{-as}

%\cnode[fillstyle=solid,fillcolor=green](1.3,3.5){0.15cm}{aa1}
  \cnode[fillstyle=solid,fillcolor=white,linewidth=1.1pt](7,0){0.12cm}{a00}
 \cnode[fillstyle=solid,fillcolor=white,linewidth=1.1pt](22,0){0.12cm}{ass}

\ncline[linewidth=1.1pt]{-}{x00}{a00}   \ncline[linewidth=1.1pt]{-}{x00}{ass}

\rput(7,-6.55){\small $0$}   \rput(14,-6.55){\small$x_0$}  \rput(22,-6.55){\small $(s-1)k$}
%\rput(9.3,4){\small $k$} \rput(12.8,4){\small $2k$}
%\rput(9,6){\small $k+1$} \rput(13.2,6){\small $2k+1$}  \rput(20.1,6.6){\small $(s-2)k+1$}

\rput(9.3,-4){\small $k$} \rput(12.8,-4){\small $2k$}  \rput(20.5,-4){\small $(s-2)k$}
\rput(8.85,-2){\small $k+1$} \rput(13.4,-2){\small $2k+1$}  \rput(19.8,-1.4){\small $(s-2)k+1$}

\rput(24.3,2){\small $(s-1)k+2$}

\ncarc[arcangle=-40,linewidth=1.3pt,linecolor=black]{-}{f1}{-a0}

\end{pspicture}
\end{center}
\caption{A drawing of the graph $H_s(k), \ k\ge 4\/$ is even.}
\end{figure}
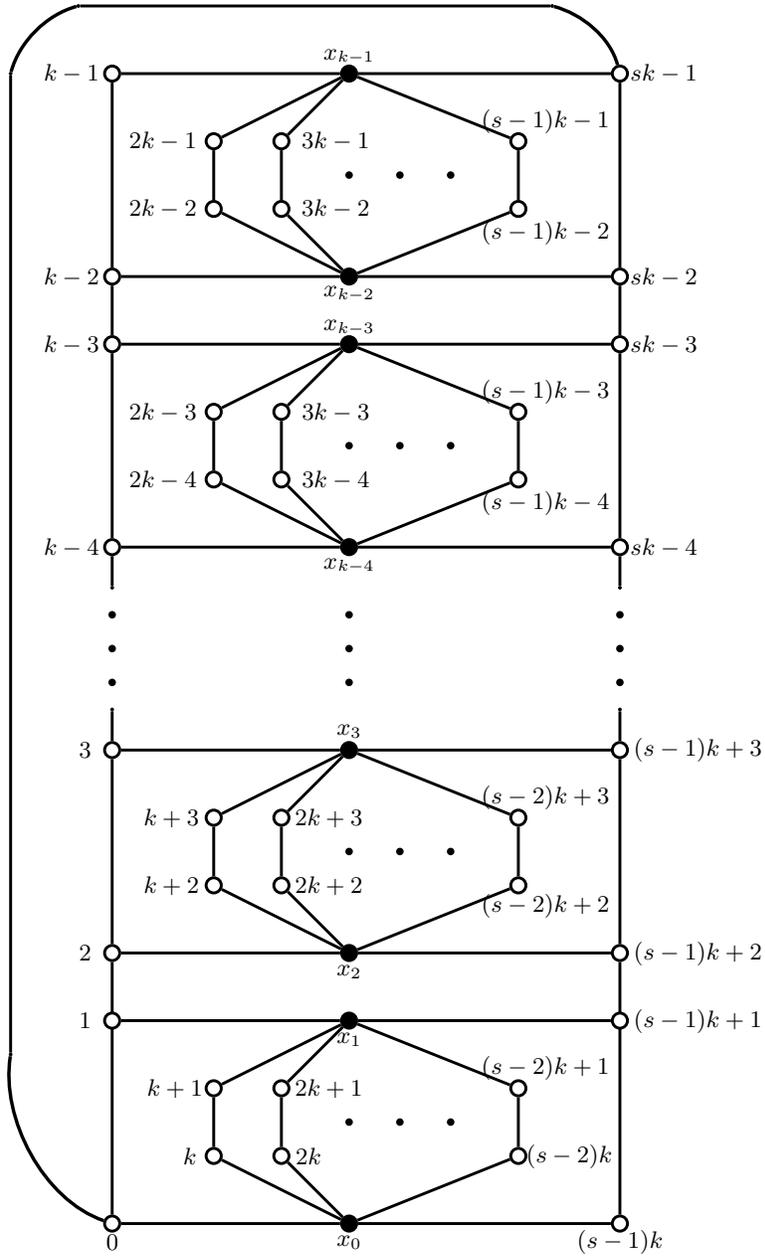

%second graph

\vspace*{8mm}
\begin{figure}[b]
\begin{center}
\psset{unit=0.36cm}
\begin{pspicture}(28,50)(3,8)
   %\showgrid

\cnode*[linewidth=1.1pt](14,6){0.12cm}{x0}     \cnode[linewidth=1.1pt](7,6){0.12cm}{c0}  \cnode[linewidth=1.1pt](22,6){0.12cm}{cs0}
\cnode*[linewidth=1.1pt](14,12){0.12cm}{x1}    \cnode[linewidth=1.1pt](7,12){0.12cm}{c1}  \cnode[linewidth=1.1pt](22,12){0.12cm}{cs1}
\cnode*[linewidth=1.1pt](14,18){0.12cm}{x2}    \cnode[linewidth=1.1pt](7,18){0.12cm}{c2}  \cnode[linewidth=1.1pt](22,18){0.12cm}{cs2}
\cnode*[linewidth=1.1pt](14,26){0.12cm}{x-4}   \cnode[linewidth=1.1pt](7,26){0.12cm}{c-4}  \cnode[linewidth=1.1pt](22,26){0.12cm}{cs-4}
\cnode*[linewidth=1.1pt](14,32){0.12cm}{x-3}   \cnode[linewidth=1.1pt](7,32){0.12cm}{c-3}  \cnode[linewidth=1.1pt](22,32){0.12cm}{cs-3}
\cnode*[linewidth=1.1pt](14,38){0.12cm}{x-2}   \cnode[linewidth=1.1pt](7,38){0.12cm}{c-2}  \cnode[linewidth=1.1pt](22,38){0.12cm}{cs-2}
\cnode*[linewidth=1.1pt](14,44){0.12cm}{x-1}   \cnode[linewidth=1.1pt](7,44){0.12cm}{c-1}  \cnode[linewidth=1.1pt](22,44){0.12cm}{cs-1}

\ncline[linewidth=1.1pt]{-}{x0}{c0}  \ncline[linewidth=1.1pt]{-}{x0}{cs0}
\ncline[linewidth=1.1pt]{-}{x1}{c1}  \ncline[linewidth=1.1pt]{-}{x1}{cs1}
\ncline[linewidth=1.1pt]{-}{x2}{c2}  \ncline[linewidth=1.1pt]{-}{x2}{cs2}
\ncline[linewidth=1.1pt]{-}{x-4}{c-4}  \ncline[linewidth=1.1pt]{-}{x-4}{cs-4}
\ncline[linewidth=1.1pt]{-}{x-3}{c-3}  \ncline[linewidth=1.1pt]{-}{x-3}{cs-3}
\ncline[linewidth=1.1pt]{-}{x-2}{c-2}  \ncline[linewidth=1.1pt]{-}{x-2}{cs-2}
\ncline[linewidth=1.1pt]{-}{x-1}{c-1}  \ncline[linewidth=1.1pt]{-}{x-1}{cs-1}

\ncline[linewidth=1.1pt]{-}{cs0}{cs1}  \ncline[linewidth=1.1pt]{-}{cs1}{cs2} \ncline[linewidth=1.1pt]{-}{cs-4}{cs-3}  \ncline[linewidth=1.1pt]{-}{cs-3}{cs-2}
\ncline[linewidth=1.1pt]{-}{cs-2}{cs-1}
\ncline[linewidth=1.1pt]{-}{c0}{c1}  \ncline[linewidth=1.1pt]{-}{c1}{c2} \ncline[linewidth=1.1pt]{-}{c-4}{c-3}  \ncline[linewidth=1.1pt]{-}{c-3}{c-2}
\ncline[linewidth=1.1pt]{-}{c-2}{c-1}

%leftc1

\cnode[linewidth=1.1pt](9,8){0.12cm}{c2k}  \cnode[linewidth=1.1pt](9,10){0.12cm}{c2k+1}
\cnode[linewidth=1.1pt](9,14){0.12cm}{ck+1}  \cnode[linewidth=1.1pt](9,16){0.12cm}{ck+2}
\cnode[linewidth=1.1pt](9,20){0.12cm}{c2k+2}  \cnode[linewidth=1.1pt](9,24){0.12cm}{c3k-4}
\cnode[linewidth=1.1pt](9,28){0.12cm}{c2k-4}  \cnode[linewidth=1.1pt](9,30){0.12cm}{c2k-3}
\cnode[linewidth=1.1pt](9,34){0.12cm}{c3k-3}  \cnode[linewidth=1.1pt](9,36){0.12cm}{c3k-2}
\cnode[linewidth=1.1pt](9,40){0.12cm}{c2k-2}  \cnode[linewidth=1.1pt](9,42){0.12cm}{c2k-1}

%leftc2

\cnode[linewidth=1.1pt](11,8){0.12cm}{c4k}  \cnode[linewidth=1.1pt](11,10){0.12cm}{c4k+1}
\cnode[linewidth=1.1pt](11,14){0.12cm}{c3k+1}  \cnode[linewidth=1.1pt](11,16){0.12cm}{c3k+2}
\cnode[linewidth=1.1pt](11,20){0.12cm}{c4k+2}  \cnode[linewidth=1.1pt](11,24){0.12cm}{c5k-4}
\cnode[linewidth=1.1pt](11,28){0.12cm}{c4k-4}  \cnode[linewidth=1.1pt](11,30){0.12cm}{c4k-3}
\cnode[linewidth=1.1pt](11,34){0.12cm}{c5k-3}  \cnode[linewidth=1.1pt](11,36){0.12cm}{c5k-2}
\cnode[linewidth=1.1pt](11,40){0.12cm}{c4k-2}  \cnode[linewidth=1.1pt](11,42){0.12cm}{c4k-1}

%rightc1

\cnode[linewidth=1.1pt](19,8){0.12cm}{cs-2k}  \cnode[linewidth=1.1pt](19,10){0.12cm}{cs-2k+1}
\cnode[linewidth=1.1pt](19,14){0.12cm}{cs-3k+1}  \cnode[linewidth=1.1pt](19,16){0.12cm}{cs-3k+2}
\cnode[linewidth=1.1pt](19,20){0.12cm}{cs-2k+2}  \cnode[linewidth=1.1pt](19,24){0.12cm}{cs-k-4}
\cnode[linewidth=1.1pt](19,28){0.12cm}{cs-2k-4}  \cnode[linewidth=1.1pt](19,30){0.12cm}{cs-2k-3}
\cnode[linewidth=1.1pt](19,34){0.12cm}{cs-k-3}  \cnode[linewidth=1.1pt](19,36){0.12cm}{cs-k-2}
\cnode[linewidth=1.1pt](19,40){0.12cm}{cs-2k-2}  \cnode[linewidth=1.1pt](19,42){0.12cm}{cs-2k-1}

%vertical line
\ncline[linewidth=1.1pt]{-}{c2k}{c2k+1} \ncline[linewidth=1.1pt]{-}{ck+1}{ck+2}  \ncline[linewidth=1.1pt]{-}{c2k-4}{c2k-3}
\ncline[linewidth=1.1pt]{-}{c3k-3}{c3k-2}  \ncline[linewidth=1.1pt]{-}{c2k-2}{c2k-1}
\ncline[linewidth=1.1pt]{-}{c4k}{c4k+1} \ncline[linewidth=1.1pt]{-}{c3k+1}{c3k+2}  \ncline[linewidth=1.1pt]{-}{c4k-4}{c4k-3}
\ncline[linewidth=1.1pt]{-}{c5k-3}{c5k-2}  \ncline[linewidth=1.1pt]{-}{c4k-2}{c4k-1}

\ncline[linewidth=1.1pt]{-}{cs-2k}{cs-2k+1} \ncline[linewidth=1.1pt]{-}{cs-3k+1}{cs-3k+2}  \ncline[linewidth=1.1pt]{-}{cs-2k-4}{cs-2k-3}
\ncline[linewidth=1.1pt]{-}{cs-k-3}{cs-k-2}  \ncline[linewidth=1.1pt]{-}{cs-2k-2}{cs-2k-1}

%diagonal
\ncline[linewidth=1.1pt]{-}{x0}{c2k} \ncline[linewidth=1.1pt]{-}{x0}{c4k} \ncline[linewidth=1.1pt]{-}{x0}{cs-2k}
\ncline[linewidth=1.1pt]{-}{x1}{c2k+1} \ncline[linewidth=1.1pt]{-}{x1}{c4k+1} \ncline[linewidth=1.1pt]{-}{x1}{cs-2k+1}
\ncline[linewidth=1.1pt]{-}{x1}{ck+1} \ncline[linewidth=1.1pt]{-}{x1}{c3k+1} \ncline[linewidth=1.1pt]{-}{x1}{cs-3k+1}
\ncline[linewidth=1.1pt]{-}{x2}{ck+2} \ncline[linewidth=1.1pt]{-}{x2}{c3k+2} \ncline[linewidth=1.1pt]{-}{x2}{cs-3k+2}
\ncline[linewidth=1.1pt]{-}{x2}{c2k+2} \ncline[linewidth=1.1pt]{-}{x2}{c4k+2} \ncline[linewidth=1.1pt]{-}{x2}{cs-2k+2}

\ncline[linewidth=1.1pt]{-}{x-4}{c3k-4} \ncline[linewidth=1.1pt]{-}{x-4}{c5k-4} \ncline[linewidth=1.1pt]{-}{x-4}{cs-k-4}
\ncline[linewidth=1.1pt]{-}{x-4}{c2k-4} \ncline[linewidth=1.1pt]{-}{x-4}{c4k-4} \ncline[linewidth=1.1pt]{-}{x-4}{cs-2k-4}
\ncline[linewidth=1.1pt]{-}{x-3}{c2k-3} \ncline[linewidth=1.1pt]{-}{x-3}{c4k-3} \ncline[linewidth=1.1pt]{-}{x-3}{cs-2k-3}
\ncline[linewidth=1.1pt]{-}{x-3}{c3k-3} \ncline[linewidth=1.1pt]{-}{x-3}{c5k-3} \ncline[linewidth=1.1pt]{-}{x-3}{cs-k-3}
\ncline[linewidth=1.1pt]{-}{x-2}{c3k-2} \ncline[linewidth=1.1pt]{-}{x-2}{c5k-2} \ncline[linewidth=1.1pt]{-}{x-2}{cs-k-2}
\ncline[linewidth=1.1pt]{-}{x-2}{c2k-2} \ncline[linewidth=1.1pt]{-}{x-2}{c4k-2} \ncline[linewidth=1.1pt]{-}{x-2}{cs-2k-2}
\ncline[linewidth=1.1pt]{-}{x-1}{c2k-1} \ncline[linewidth=1.1pt]{-}{x-1}{c4k-1} \ncline[linewidth=1.1pt]{-}{x-1}{cs-2k-1}

%small node
\cnode*[linewidth=1.1pt](7,19.5){0.025cm}{h1}    \cnode*[linewidth=1.1pt](7,24.5){0.025cm}{h2}
\cnode*[linewidth=1.1pt](9,20.5){0.025cm}{h3}    \cnode*[linewidth=1.1pt](9,23.5){0.025cm}{h4}
\cnode*[linewidth=1.1pt](11,20.5){0.025cm}{h5}   \cnode*[linewidth=1.1pt](11,23.5){0.025cm}{h6}
\cnode*[linewidth=1.1pt](19,20.5){0.025cm}{h7}   \cnode*[linewidth=1.1pt](19,23.5){0.025cm}{h8}
\cnode*[linewidth=1.1pt](22,19.5){0.025cm}{h9}    \cnode*[linewidth=1.1pt](22,24.5){0.025cm}{h10}

\ncline[linewidth=1.1pt]{-}{h2}{c-4} \ncline[linewidth=1.1pt]{-}{h1}{c2} \ncline[linewidth=1.1pt]{-}{h10}{cs-4} \ncline[linewidth=1.1pt]{-}{h9}{cs2}
\ncline[linewidth=1.1pt]{-}{h4}{c3k-4} \ncline[linewidth=1.1pt]{-}{h3}{c2k+2} \ncline[linewidth=1.1pt]{-}{h5}{c4k+2} \ncline[linewidth=1.1pt]{-}{h6}{c5k-4}
\ncline[linewidth=1.1pt]{-}{h8}{cs-k-4} \ncline[linewidth=1.1pt]{-}{h7}{cs-2k+2}

\cnode*[linewidth=1.1pt](7,21){0.05cm}{d0} \cnode*[linewidth=1.1pt](7,22){0.05cm}{d1} \cnode*[linewidth=1.1pt](7,23){0.05cm}{d2}
\cnode*[linewidth=1.1pt](22,21){0.05cm}{d0} \cnode*[linewidth=1.1pt](22,22){0.05cm}{d1} \cnode*[linewidth=1.1pt](22,23){0.05cm}{d2}

\cnode*[linewidth=1.1pt](9,21){0.05cm}{d0} \cnode*[linewidth=1.1pt](9,22){0.05cm}{d1} \cnode*[linewidth=1.1pt](9,23){0.05cm}{d2}
\cnode*[linewidth=1.1pt](11,21){0.05cm}{d0} \cnode*[linewidth=1.1pt](11,22){0.05cm}{d1} \cnode*[linewidth=1.1pt](11,23){0.05cm}{d2}
\cnode*[linewidth=1.1pt](19,21){0.05cm}{d0} \cnode*[linewidth=1.1pt](19,22){0.05cm}{d1} \cnode*[linewidth=1.1pt](19,23){0.05cm}{d2}
\cnode*[linewidth=1.1pt](3,21){0.05cm}{d0} \cnode*[linewidth=1.1pt](3,22){0.05cm}{d1} \cnode*[linewidth=1.1pt](3,23){0.05cm}{d2}
\cnode*[linewidth=1.1pt](4,21){0.05cm}{d0} \cnode*[linewidth=1.1pt](4,22){0.05cm}{d1} \cnode*[linewidth=1.1pt](4,23){0.05cm}{d2}
\cnode*[linewidth=1.1pt](25,21){0.05cm}{d0} \cnode*[linewidth=1.1pt](25,22){0.05cm}{d1} \cnode*[linewidth=1.1pt](25,23){0.05cm}{d2}
\cnode*[linewidth=1.1pt](26,21){0.05cm}{d0} \cnode*[linewidth=1.1pt](26,22){0.05cm}{d1} \cnode*[linewidth=1.1pt](26,23){0.05cm}{d2}

\cnode[linewidth=1.1pt](7,4){0.12cm}{ck}  \cnode[linewidth=1.1pt](8,2){0.12cm}{c3k} \cnode[linewidth=1.1pt](22,4){0.12cm}{cs-3k}
\ncline[linewidth=1.1pt]{-}{x0}{ck} \ncline[linewidth=1.1pt]{-}{x0}{c3k} \ncline[linewidth=1.1pt]{-}{x0}{cs-3k}

\cnode[linewidth=1.1pt](9,46){0.12cm}{c3k-1}  \cnode[linewidth=1.1pt](22,46){0.12cm}{cs-k-1} \cnode[linewidth=1.1pt](21,47.5){0.12cm}{cs-3k-1}
\ncline[linewidth=1.1pt]{-}{x-1}{c3k-1} \ncline[linewidth=1.1pt]{-}{x-1}{cs-3k-1} \ncline[linewidth=1.1pt]{-}{x-1}{cs-k-1}

\cnode*[linewidth=1.1pt](4,7){0.025cm}{l1}  \cnode*[linewidth=1.1pt](4,19.5){0.025cm}{l2}  \cnode*[linewidth=1.1pt](4,24.5){0.025cm}{l3} \cnode*[linewidth=1.1pt](4,41){0.025cm}{l4}
\cnode*[linewidth=1.1pt](3,7){0.025cm}{l5}  \cnode*[linewidth=1.1pt](3,19.5){0.025cm}{l6}  \cnode*[linewidth=1.1pt](3,24.5){0.025cm}{l7} \cnode*[linewidth=1.1pt](3,41){0.025cm}{l8}

\cnode*[linewidth=1.1pt](8.5,46){0.025cm}{l9}  \cnode*[linewidth=1.1pt](7,2.2){0.025cm}{l10}  \cnode*[linewidth=1.1pt](25,43){0.025cm}{l11} \cnode*[linewidth=1.1pt](25,9){0.025cm}{l12} \cnode*[linewidth=1.1pt](25,46){0.025cm}{l13}  \cnode*[linewidth=1.1pt](26,43){0.025cm}{l14}
\cnode*[linewidth=1.1pt](26,9){0.025cm}{l15} \cnode*[linewidth=1.1pt](25,24.5){0.025cm}{l16}  \cnode*[linewidth=1.1pt](25,19.5){0.025cm}{l17}
\cnode*[linewidth=1.1pt](26,24.5){0.025cm}{l18}  \cnode*[linewidth=1.1pt](26,19.5){0.025cm}{l19}

\ncline[linewidth=1.1pt]{-}{l1}{l2}  \ncline[linewidth=1.1pt]{-}{l3}{l4} \ncline[linewidth=1.1pt]{-}{l5}{l6}  \ncline[linewidth=1.1pt]{-}{l7}{l8}
\ncline[linewidth=1.1pt]{-}{l11}{l16}  \ncline[linewidth=1.1pt]{-}{l17}{l12} \ncline[linewidth=1.1pt]{-}{l14}{l18}  \ncline[linewidth=1.1pt]{-}{l15}{l19}

\ncarc[arcangle=30,linewidth=1.3pt,linecolor=black]{-}{l4}{c-1}  \ncarc[arcangle=-30,linewidth=1.3pt,linecolor=black]{-}{l1}{ck}  \ncarc[arcangle=-30,linewidth=1.3pt,linecolor=black]{-}{l11}{cs-k-1}  \ncarc[arcangle=-30,linewidth=1.3pt,linecolor=black]{-}{l13}{cs-3k-1}
\ncarc[arcangle=15,linewidth=1.3pt,linecolor=black]{-}{l13}{l14}  \ncarc[arcangle=-30,linewidth=1.3pt,linecolor=black]{-}{cs0}{l12}
\ncarc[arcangle=-30,linewidth=1.3pt,linecolor=black]{-}{cs-3k}{l15}

\ncline[linewidth=1.1pt]{-}{l9}{c3k-1} \ncarc[arcangle=-30,linewidth=1.3pt,linecolor=black]{-}{l9}{l8}
\ncline[linewidth=1.1pt]{-}{l10}{c3k} \ncarc[arcangle=30,linewidth=1.3pt,linecolor=black]{-}{l10}{l5}

\cnode*[linewidth=1.1pt](14,46.2){0.05cm}{d0} \cnode*[linewidth=1.1pt](12.5,46){0.05cm}{d1} \cnode*[linewidth=1.1pt](15.5,46){0.05cm}{d2}
\cnode*[linewidth=1.1pt](14.5,3.8){0.05cm}{d0} \cnode*[linewidth=1.1pt](16,4){0.05cm}{d1} \cnode*[linewidth=1.1pt](13,4){0.05cm}{d2}

\rput(14,6.5){\small $x_0$} \rput(14,12.5){\small $x_1$} \rput(14,18.5){\small $x_2$} \rput(14,26.7){\small $x_{k-4}$} \rput(14,32.7){\small $x_{k-3}$}
\rput(14,38.7){\small $x_{k-2}$} \rput(14,44.7){\small $x_{k-1}$}

\rput(6.3,6){\small $0$} \rput(6.3,12){\small $1$} \rput(6.3,18){\small $2$} \rput(5.65,26){\small $k-4$} \rput(5.65,32){\small $k-3$} \rput(5.65,38){\small $k-2$}
\rput(7,44.5){\small $k-1$} \rput(7,3.4){\small $k$}  \rput(8,2.75){\small $3k$}

\rput(8,8){\small $2k$} \rput(8.3,10.6){\small $2k+1$} \rput(8.3,13.4){\small $k+1$} \rput(8.5,16.6){\small $k+2$} \rput(8.3,19.4){\small $2k+2$} \rput(8.3,24.6){\small $3k-4$}
\rput(8.3,27.5){\small $2k-4$} \rput(8.3,30.5){\small $2k-3$} \rput(8.5,33.4){\small $3k-3$} \rput(8.3,36.6){\small $3k-2$} \rput(8.3,39.4){\small $2k-2$}
\rput(8.3,42.6){\small $2k-1$} \rput(9,46.5){\small $3k-1$}

\rput(12.1,8){\small $4k$}   \rput(12.7,10){\small $4k+1$}  \rput(12.7,14){\small $3k+1$} \rput(12.7,16){\small $3k+2$} \rput(12.7,20){\small $4k+2$} \rput(12.7,24){\small $5k-4$}
\rput(12.7,28){\small $4k-4$}  \rput(12.7,30){\small $4k-3$} \rput(12.7,34){\small $5k-3$} \rput(12.8,35.8){\small $5k-2$} \rput(12.7,40){\small $4k-2$} \rput(12.7,42){\small $4k-1$}

\rput(22,3.4){\small $(s-3)k$}  \rput(22,5.4){\small $(s-1)k$} \rput(24.3,12.3){\small $(s-1)k+1$} \rput(24.3,18.3){\small $(s-1)k+2$}   \rput(23.5,26){\small $sk-4$}
 \rput(23.5,32){\small $sk-3$}    \rput(23.5,38){\small $sk-2$}  \rput(22,44.5){\small $sk-1$}

 \rput(20,7.4){\small $(s-2)k$}  \rput(19.5,10.8){\small $(s-2)k+1$} \rput(19.5,13.2){\small $(s-3)k+1$} \rput(19.5,16.8){\small $(s-3)k+2$}
\rput(19.5,19.2){\small $(s-2)k+2$} \rput(19.5,24.8){\small $(s-1)k-4$}  \rput(19.5,27.2){\small $(s-2)k-4$} \rput(19.5,30.8){\small $(s-2)k-3$}
\rput(19.5,33.2){\small $(s-1)k-3$} \rput(19.5,36.8){\small $(s-1)k-2$}  \rput(19.6,39.1){\small $(s-2)k-2$} \rput(19.6,42.8){\small $(s-2)k-1$}

\rput(18.2,47.5){\small $(s-3)k-1$}  \rput(22,46.5){\small $(s-1)k-1$}

\cnode*[linewidth=1.1pt](14,15){0.05cm}{d9} \cnode*[linewidth=1.1pt](15.5,15){0.05cm}{d10} \cnode*[linewidth=1.1pt](17,15){0.05cm}{d11}
\cnode*[linewidth=1.1pt](14,9){0.05cm}{d9} \cnode*[linewidth=1.1pt](15.5,9){0.05cm}{d10} \cnode*[linewidth=1.1pt](17,9){0.05cm}{d11}
\cnode*[linewidth=1.1pt](14,29){0.05cm}{d9} \cnode*[linewidth=1.1pt](15.5,29){0.05cm}{d10} \cnode*[linewidth=1.1pt](17,29){0.05cm}{d11}
\cnode*[linewidth=1.1pt](14,35){0.05cm}{d9} \cnode*[linewidth=1.1pt](15.5,35){0.05cm}{d10} \cnode*[linewidth=1.1pt](17,35){0.05cm}{d11}
\cnode*[linewidth=1.1pt](14,41){0.05cm}{d9} \cnode*[linewidth=1.1pt](15.5,41){0.05cm}{d10} \cnode*[linewidth=1.1pt](17,41){0.05cm}{d11}
\cnode*[linewidth=1.1pt](14,20){0.05cm}{d9} \cnode*[linewidth=1.1pt](15.5,20){0.05cm}{d10} \cnode*[linewidth=1.1pt](17,20){0.05cm}{d11}
\cnode*[linewidth=1.1pt](14,24){0.05cm}{d9} \cnode*[linewidth=1.1pt](15.5,24){0.05cm}{d10} \cnode*[linewidth=1.1pt](17,24){0.05cm}{d11}

\end{pspicture}
\end{center}

\vspace{18mm}

\caption{A drawing of the graph $H_s(k), \ k\ge 5\/$ is odd and $s \ge 4\/$ is even.}
\end{figure}
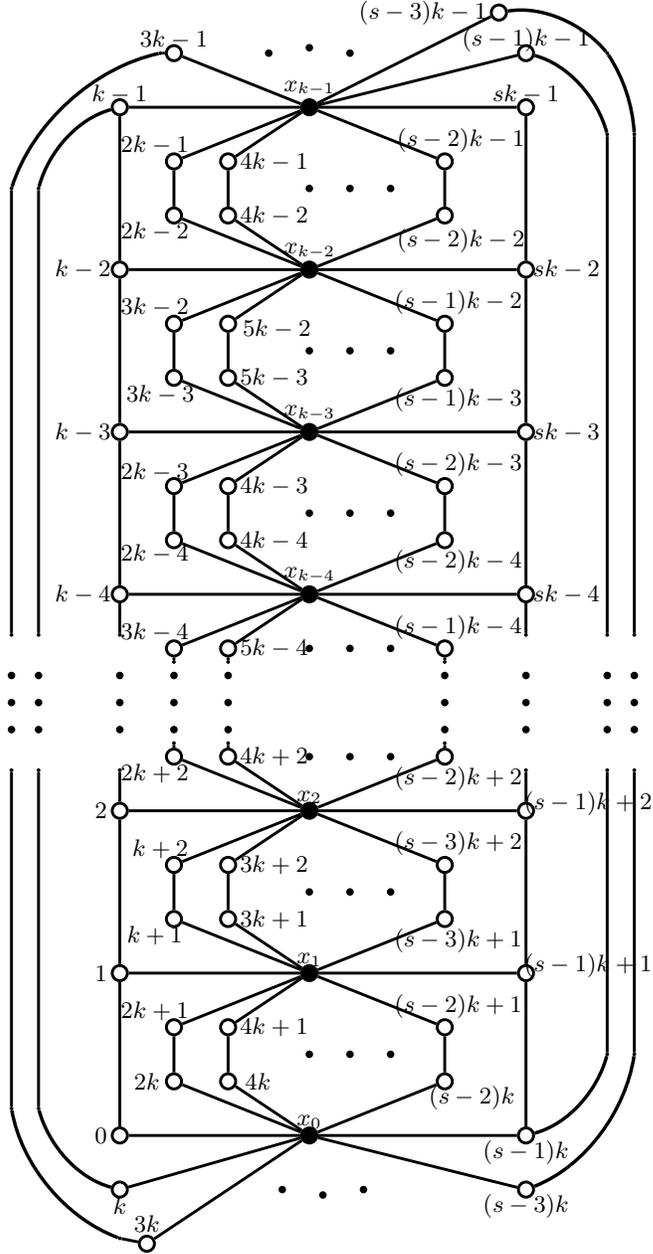

%third graph

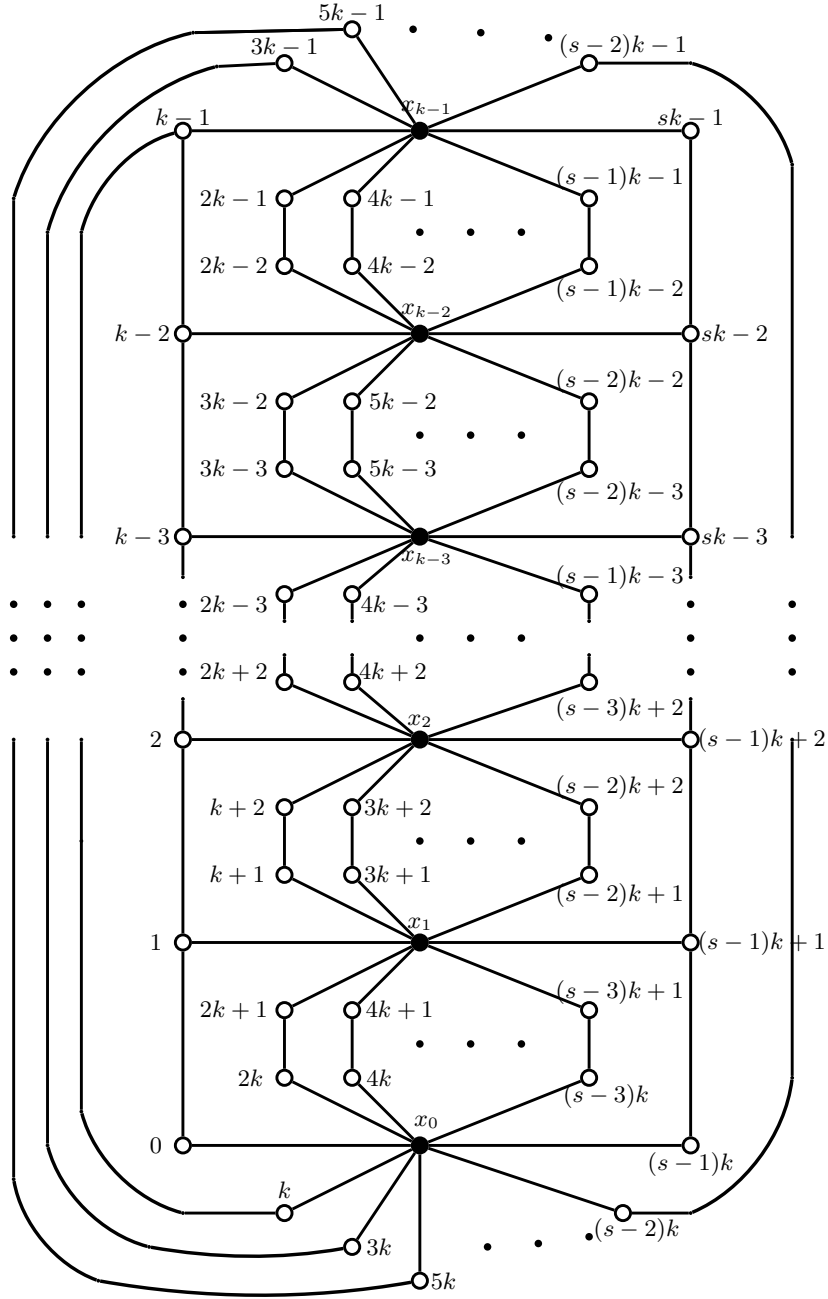
\begin{figure}
\begin{center}
\psset{unit=0.45cm}
\begin{pspicture}(28,30)(1,-8)
   %\showgrid

\cnode[linewidth=1.1pt](7,-4){0.12cm}{-a0} \cnode*[linewidth=1.1pt](14,-4){0.12cm}{-x0}  \cnode[linewidth=1.1pt](22,-4){0.12cm}{-as}
\ncline[linewidth=1.1pt]{-}{-x0}{-a0}  \ncline[linewidth=1.1pt]{-}{-x0}{-as}

\cnode*[linewidth=1.1pt](14,2){0.12cm}{x0}
\cnode*[linewidth=1.1pt](14,8){0.12cm}{x1}
\cnode*[linewidth=1.1pt](14,14){0.12cm}{x-4}
\cnode*[linewidth=1.1pt](14,20){0.12cm}{x-3}
\cnode*[linewidth=1.1pt](14,26){0.12cm}{x-2}
%\cnode*[linewidth=1.1pt](14,32){0.12cm}{x-1}
%\cnode*[linewidth=1.1pt](14,38){0.12cm}{x}

\cnode[linewidth=1.1pt](7,2){0.12cm}{a0}
\cnode[linewidth=1.1pt](7,8){0.12cm}{a1}
\cnode[linewidth=1.1pt](7,14){0.12cm}{a-4}
\cnode[linewidth=1.1pt](7,20){0.12cm}{a-3}
\cnode[linewidth=1.1pt](7,26){0.12cm}{a-2}
%\cnode[linewidth=1.1pt](7,32){0.12cm}{a-1}
%\cnode[linewidth=1.1pt](7,38){0.12cm}{a}

\cnode[linewidth=1.1pt](22,2){0.12cm}{as}
\cnode[linewidth=1.1pt](22,8){0.12cm}{as1}
\cnode[linewidth=1.1pt](22,14){0.12cm}{as-4}
\cnode[linewidth=1.1pt](22,20){0.12cm}{as-3}
\cnode[linewidth=1.1pt](22,26){0.12cm}{as-2}
%\cnode[linewidth=1.1pt](22,32){0.12cm}{as-1}
%\cnode[linewidth=1.1pt](22,38){0.12cm}{as01}

\ncline[linewidth=1.1pt]{-}{a0}{-a0}  \ncline[linewidth=1.1pt]{-}{as}{-as}

\cnode*[linewidth=1.1pt](7,10){0.05cm}{d0} \cnode*[linewidth=1.1pt](7,11){0.05cm}{d1} \cnode*[linewidth=1.1pt](7,12){0.05cm}{d2}
\cnode*[linewidth=1.1pt](22,10){0.05cm}{d3} \cnode*[linewidth=1.1pt](22,11){0.05cm}{d4} \cnode*[linewidth=1.1pt](22,12){0.05cm}{d5}
%\cnode*[linewidth=1.1pt](14,-1){0.05cm}{d6} \cnode*[linewidth=1.1pt](14,-1){0.05cm}{d7} \cnode*[linewidth=1.1pt](14,12){0.05cm}{d8}

\cnode*[linewidth=1.1pt](14,17){0.05cm}{d9} \cnode*[linewidth=1.1pt](15.5,17){0.05cm}{d10} \cnode*[linewidth=1.1pt](17,17){0.05cm}{d11}
\cnode*[linewidth=1.1pt](14,5){0.05cm}{d9} \cnode*[linewidth=1.1pt](15.5,5){0.05cm}{d10} \cnode*[linewidth=1.1pt](17,5){0.05cm}{d11}
\cnode*[linewidth=1.1pt](14,23){0.05cm}{d9} \cnode*[linewidth=1.1pt](15.5,23){0.05cm}{d10} \cnode*[linewidth=1.1pt](17,23){0.05cm}{d11}
\cnode*[linewidth=1.1pt](14,11){0.05cm}{d12} \cnode*[linewidth=1.1pt](15.5,11){0.05cm}{d13} \cnode*[linewidth=1.1pt](17,11){0.05cm}{d14}
\cnode*[linewidth=1.1pt](14,-1){0.05cm}{d12} \cnode*[linewidth=1.1pt](15.5,-1){0.05cm}{d13} \cnode*[linewidth=1.1pt](17,-1){0.05cm}{d14}

\cnode*[linewidth=1.1pt](7,9.2){0.025cm}{e1} \cnode*[linewidth=1.1pt](7,12.8){0.025cm}{e2}
\cnode*[linewidth=1.1pt](22,9.2){0.025cm}{e3} \cnode*[linewidth=1.1pt](22,12.8){0.025cm}{e4}
\ncline[linewidth=1.1pt]{-}{a1}{e1} \ncline[linewidth=1.1pt]{-}{a-4}{e2}  \ncline[linewidth=1.1pt]{-}{as1}{e3} \ncline[linewidth=1.1pt]{-}{as-4}{e4}

\cnode*[linewidth=1.1pt](4,5){0.025cm}{f1}

%\ncline[linewidth=1.1pt]{-}{f1}{f2}  \ncline[linewidth=1.1pt]{-}{f3}{f4}
%\ncarc[arcangle=-40,linewidth=1.3pt,linecolor=black]{-}{f1}{a0} \ncarc[arcangle=30,linewidth=1.3pt,linecolor=black]{-}{f2}{f3} \ncarc[arcangle=30,linewidth=1.3pt,linecolor=black]{-}{f4}{as-1}

\ncline[linewidth=1.1pt]{-}{x0}{a0}  \ncline[linewidth=1.1pt]{-}{x0}{as}
\ncline[linewidth=1.1pt]{-}{x1}{a1}  \ncline[linewidth=1.1pt]{-}{x1}{as1}
\ncline[linewidth=1.1pt]{-}{x-4}{a-4}  \ncline[linewidth=1.1pt]{-}{x-4}{as-4}
\ncline[linewidth=1.1pt]{-}{x-3}{a-3}  \ncline[linewidth=1.1pt]{-}{x-3}{as-3}
\ncline[linewidth=1.1pt]{-}{x-2}{a-2}  \ncline[linewidth=1.1pt]{-}{x-2}{as-2}
\ncline[linewidth=1.1pt]{-}{x-1}{a-1}  \ncline[linewidth=1.1pt]{-}{x-1}{as-1}

\ncline[linewidth=1.1pt]{-}{a0}{a1}  \ncline[linewidth=1.1pt]{-}{a-4}{a-3}  \ncline[linewidth=1.1pt]{-}{a-3}{a-2}  \ncline[linewidth=1.1pt]{-}{a-2}{a-1}
\ncline[linewidth=1.1pt]{-}{as}{as1}  \ncline[linewidth=1.1pt]{-}{as-4}{as-3}  \ncline[linewidth=1.1pt]{-}{as-3}{as-2}  \ncline[linewidth=1.1pt]{-}{as-2}{as-1}

\cnode[linewidth=1.1pt](10,4){0.12cm}{ak} \cnode[linewidth=1.1pt](12,4){0.12cm}{a2k}   \cnode[linewidth=1.1pt](19,4){0.12cm}{as-2k}
\cnode[linewidth=1.1pt](10,6){0.12cm}{ak1} \cnode[linewidth=1.1pt](12,6){0.12cm}{a2k1}   \cnode[linewidth=1.1pt](19,6){0.12cm}{as-2k1}

\ncline[linewidth=1.1pt]{-}{x0}{ak}   \ncline[linewidth=1.1pt]{-}{x0}{a2k} \ncline[linewidth=1.1pt]{-}{x0}{as-2k}
\ncline[linewidth=1.1pt]{-}{x1}{ak1}   \ncline[linewidth=1.1pt]{-}{x1}{a2k1} \ncline[linewidth=1.1pt]{-}{x1}{as-2k1}
\ncline[linewidth=1.1pt]{-}{ak1}{ak}   \ncline[linewidth=1.1pt]{-}{a2k1}{a2k} \ncline[linewidth=1.1pt]{-}{as-2k1}{as-2k}

\cnode[linewidth=1.1pt](10,16){0.12cm}{a2k-4} \cnode[linewidth=1.1pt](12,16){0.12cm}{a3k-4}   \cnode[linewidth=1.1pt](19,16){0.12cm}{ask-4}
\cnode[linewidth=1.1pt](10,18){0.12cm}{a2k-3} \cnode[linewidth=1.1pt](12,18){0.12cm}{a3k-3}   \cnode[linewidth=1.1pt](19,18){0.12cm}{ask-3}

\ncline[linewidth=1.1pt]{-}{x-4}{a2k-4}   \ncline[linewidth=1.1pt]{-}{x-4}{a3k-4} \ncline[linewidth=1.1pt]{-}{x-4}{ask-4}
\ncline[linewidth=1.1pt]{-}{x-3}{a2k-3}   \ncline[linewidth=1.1pt]{-}{x-3}{a3k-3} \ncline[linewidth=1.1pt]{-}{x-3}{ask-3}
\ncline[linewidth=1.1pt]{-}{a2k-3}{a2k-4}   \ncline[linewidth=1.1pt]{-}{a3k-3}{a3k-4} \ncline[linewidth=1.1pt]{-}{ask-3}{ask-4}

\cnode[linewidth=1.1pt](10,28){0.12cm}{a2k-2} \cnode[linewidth=1.1pt](12,29){0.12cm}{a3k-2}   \cnode[linewidth=1.1pt](19,28){0.12cm}{ask-2}
%\cnode[linewidth=1.1pt](10,30){0.12cm}{a2k-1} \cnode[linewidth=1.1pt](12,30){0.12cm}{a3k-1}   \cnode[linewidth=1.1pt](19,30){0.12cm}{ask-1}

\ncline[linewidth=1.1pt]{-}{x-2}{a2k-2}   \ncline[linewidth=1.1pt]{-}{x-2}{a3k-2} \ncline[linewidth=1.1pt]{-}{x-2}{ask-2}
\ncline[linewidth=1.1pt]{-}{x-1}{a2k-1}   \ncline[linewidth=1.1pt]{-}{x-1}{a3k-1} \ncline[linewidth=1.1pt]{-}{x-1}{ask-1}
\ncline[linewidth=1.1pt]{-}{a2k-2}{a2k-1}   \ncline[linewidth=1.1pt]{-}{a3k-2}{a3k-1} \ncline[linewidth=1.1pt]{-}{ask-2}{ask-1}

%\rput(14,-10){(a) $k\/$ and $s\/$ are odd.}

\rput(6.2,-4){\small $0$}  \rput(6.2,2){\small $1$}  \rput(14.2,-3.3){\small$x_0$}  \rput(22,-4.55){\small $(s-1)k$}
\rput(6.2,8){\small $2$}   \rput(14,2.55){\small$x_1$}  \rput(24.1,2){\small $(s-1)k+1$}  \rput(14,8.55){\small$x_2$}
\rput(5.8,14){\small $k-3$}   \rput(14.2,13.35){\small$x_{k-3}$}  \rput(23.3,14){\small $sk-3$}
\rput(5.8,20){\small $k-2$}   \rput(14.2,20.65){\small$x_{k-2}$}  %\rput(23.3,20){\small $sk-2$}
    \rput(23.3,20){\small $sk-2$} %\rput(5.8,22){\small $k-2$} \rput(14,21.5){\small$x_{k-2}$}
\rput(7,26.4){\small $k-1$}   \rput(14.2,26.65){\small$x_{k-1}$}  \rput(22,26.4){\small $sk-1$}
\rput(19.5,-2.55){\small $(s-3)k$}  \rput(19.9,0.55){\small $(s-3)k+1$}   \rput(24.1,8){\small $(s-1)k+2$}

\rput(8.6,4){\small $k+1$} \rput(13.35,4){\small $3k+1$}  \rput(19.9,3.4){\small $(s-2)k+1$}
\rput(8.6,6){\small $k+2$} \rput(13.35,6){\small $3k+2$}  \rput(19.9,6.6){\small $(s-2)k+2$}
\rput(19.9,12.8){\small $(s-1)k-3$}  \rput(19.9,8.9){\small $(s-3)k+2$}

\rput(8.5,16){\small $3k-3$} \rput(13.5,16){\small $5k-3$}  \rput(19.9,15.3){\small $(s-2)k-3$}
\rput(8.5,18){\small $3k-2$} \rput(13.5,18){\small $5k-2$}  \rput(19.9,18.6){\small $(s-2)k-2$}
\rput(8.5,12){\small $2k-3$}  \rput(8.5,10){\small $2k+2$}
\rput(8.5,0){\small $2k+1$}  \rput(9,-2){\small $2k$}
\rput(13.2,10){\small $4k+2$} \rput(13.25,12){\small $4k-3$}  \rput(13.4,0){\small $4k+1$} \rput(12.8,-2){\small $4k$}

\rput(8.5,24){\small $2k-1$} \rput(13.45,24){\small $4k-1$}  \rput(19.9,21.3){\small $(s-1)k-2$}
\rput(8.5,22){\small $2k-2$} \rput(13.45,22){\small $4k-2$}  \rput(19.9,24.6){\small $(s-1)k-1$}

\ncline[linewidth=1.1pt]{-}{x}{a} \ncline[linewidth=1.1pt]{-}{x}{as01}  \ncline[linewidth=1.1pt]{-}{a-1}{a} \ncline[linewidth=1.1pt]{-}{as-1}{as01}

%\cnode[linewidth=1.1pt](10,34){0.12cm}{aa2k-2} \cnode[linewidth=1.1pt](12,34){0.12cm}{aa3k-2}   \cnode[linewidth=1.1pt](19,34){0.12cm}{aask-2}
%\cnode[linewidth=1.1pt](10,36){0.12cm}{aa2k-1} \cnode[linewidth=1.1pt](12,36){0.12cm}{aa3k-1}   \cnode[linewidth=1.1pt](19,36){0.12cm}{aask-1}

\ncline[linewidth=1.1pt]{-}{x}{aa2k-1} \ncline[linewidth=1.1pt]{-}{x}{aa3k-1}   \ncline[linewidth=1.1pt]{-}{x}{aask-1}
\ncline[linewidth=1.1pt]{-}{x-1}{aa2k-2} \ncline[linewidth=1.1pt]{-}{x-1}{aa3k-2}   \ncline[linewidth=1.1pt]{-}{x-1}{aask-2}
\ncline[linewidth=1.1pt]{-}{aa2k-2}{aa2k-1} \ncline[linewidth=1.1pt]{-}{aa3k-2}{aa3k-1}   \ncline[linewidth=1.1pt]{-}{aask-2}{aask-1}

\cnode[linewidth=1.1pt](10,22){0.12cm}{ab2k-2} \cnode[linewidth=1.1pt](12,22){0.12cm}{ab3k-2}  \cnode[linewidth=1.1pt](19,22){0.12cm}{absk-2}
\cnode[linewidth=1.1pt](10,24){0.12cm}{ab2k-1} \cnode[linewidth=1.1pt](12,24){0.12cm}{ab3k-1}  \cnode[linewidth=1.1pt](19,24){0.12cm}{absk-1}

\ncline[linewidth=1.1pt]{-}{x-2}{ab2k-1} \ncline[linewidth=1.1pt]{-}{x-2}{ab3k-1}  \ncline[linewidth=1.1pt]{-}{x-2}{absk-1}
\ncline[linewidth=1.1pt]{-}{x-3}{ab2k-2} \ncline[linewidth=1.1pt]{-}{x-3}{ab3k-2}  \ncline[linewidth=1.1pt]{-}{x-3}{absk-2}
\ncline[linewidth=1.1pt]{-}{ab2k-1}{ab2k-2} \ncline[linewidth=1.1pt]{-}{ab3k-1}{ab3k-2}  \ncline[linewidth=1.1pt]{-}{absk-1}{absk-2}

\cnode[linewidth=1.1pt](10,0){0.12cm}{ac1} \cnode[linewidth=1.1pt](12,0){0.12cm}{ac2}  \cnode[linewidth=1.1pt](19,0){0.12cm}{ac3}
\cnode[linewidth=1.1pt](10,-2){0.12cm}{ac-1} \cnode[linewidth=1.1pt](12,-2){0.12cm}{ac-2}  \cnode[linewidth=1.1pt](19,-2){0.12cm}{ac-3}

\ncline[linewidth=1.1pt]{-}{x0}{ac1} \ncline[linewidth=1.1pt]{-}{x0}{ac2}  \ncline[linewidth=1.1pt]{-}{x0}{ac3}
\ncline[linewidth=1.1pt]{-}{-x0}{ac-1} \ncline[linewidth=1.1pt]{-}{-x0}{ac-2}  \ncline[linewidth=1.1pt]{-}{-x0}{ac-3}
\ncline[linewidth=1.1pt]{-}{ac1}{ac-1} \ncline[linewidth=1.1pt]{-}{ac2}{ac-2}  \ncline[linewidth=1.1pt]{-}{ac3}{ac-3}

\cnode[linewidth=1.1pt](10,12.3){0.12cm}{ca1} \cnode[linewidth=1.1pt](12,12.3){0.12cm}{ca2}  \cnode[linewidth=1.1pt](19,12.3){0.12cm}{ca3}
\cnode[linewidth=1.1pt](10,9.7){0.12cm}{ca-1} \cnode[linewidth=1.1pt](12,9.7){0.12cm}{ca-2}  \cnode[linewidth=1.1pt](19,9.7){0.12cm}{ca-3}

\ncline[linewidth=1.1pt]{-}{x1}{ca-1} \ncline[linewidth=1.1pt]{-}{x1}{ca-2}  \ncline[linewidth=1.1pt]{-}{x1}{ca-3}
\ncline[linewidth=1.1pt]{-}{x-4}{ca1} \ncline[linewidth=1.1pt]{-}{x-4}{ca2}  \ncline[linewidth=1.1pt]{-}{x-4}{ca3}

\cnode*[linewidth=1.1pt](10,11.5){0.025cm}{f2} \cnode*[linewidth=1.1pt](10,10.5){0.025cm}{f3}
\cnode*[linewidth=1.1pt](19,11.5){0.025cm}{f4} \cnode*[linewidth=1.1pt](19,10.5){0.025cm}{f5}
\cnode*[linewidth=1.1pt](12,11.5){0.025cm}{f6} \cnode*[linewidth=1.1pt](12,10.5){0.025cm}{f7}

\ncline[linewidth=1.1pt]{-}{f2}{ca1}  \ncline[linewidth=1.1pt]{-}{f3}{ca-1}  \ncline[linewidth=1.1pt]{-}{f4}{ca3}  \ncline[linewidth=1.1pt]{-}{f5}{ca-3}
\ncline[linewidth=1.1pt]{-}{f6}{ca2}  \ncline[linewidth=1.1pt]{-}{f7}{ca-2}

\cnode*[linewidth=1.1pt](13.8,29){0.05cm}{d9}  \cnode*[linewidth=1.1pt](15.8,28.9){0.05cm}{d9}  \cnode*[linewidth=1.1pt](17.8,28.75){0.05cm}{d9}

\cnode[linewidth=1.1pt](10,-6){0.12cm}{k1}   \cnode[linewidth=1.1pt](20,-6){0.12cm}{ks}   \cnode[linewidth=1.1pt](12,-7){0.12cm}{k3}   \cnode[linewidth=1.1pt](14,-8){0.12cm}{k5}
\cnode*[linewidth=1.1pt](16,-7){0.05cm}{d9}  \cnode*[linewidth=1.1pt](19,-6.7){0.05cm}{d9}  \cnode*[linewidth=1.1pt](17.5,-6.9){0.05cm}{d9}
\ncline[linewidth=1.1pt]{-}{-x0}{k1} \ncline[linewidth=1.1pt]{-}{-x0}{k3}   \ncline[linewidth=1.1pt]{-}{-x0}{k5}  \ncline[linewidth=1.1pt]{-}{-x0}{ks}

\cnode*[linewidth=1.1pt](2,-5){0.025cm}{g1} \cnode*[linewidth=1.1pt](3,-4){0.025cm}{g2}  \cnode*[linewidth=1.1pt](4,-3){0.025cm}{g3} \cnode*[linewidth=1.1pt](25,-2){0.025cm}{g4}
\cnode*[linewidth=1.1pt](2,8){0.025cm}{g5} \cnode*[linewidth=1.1pt](3,8){0.025cm}{g6}  \cnode*[linewidth=1.1pt](4,8){0.025cm}{g7} \cnode*[linewidth=1.1pt](25,8){0.025cm}{g8}
\cnode*[linewidth=1.1pt](2,14){0.025cm}{g9} \cnode*[linewidth=1.1pt](3,14){0.025cm}{g10}  \cnode*[linewidth=1.1pt](4,14){0.025cm}{g11} \cnode*[linewidth=1.1pt](25,14){0.025cm}{g12}
\cnode*[linewidth=1.1pt](2,24){0.025cm}{g13} \cnode*[linewidth=1.1pt](3,23){0.025cm}{g14}  \cnode*[linewidth=1.1pt](4,23){0.025cm}{g15} \cnode*[linewidth=1.1pt](25,25){0.025cm}{g16}
\cnode*[linewidth=1.1pt](7.3,28.9){0.025cm}{g17} \cnode*[linewidth=1.1pt](8,27.9){0.025cm}{g18}  \cnode*[linewidth=1.1pt](22,28){0.025cm}{g19}

\cnode*[linewidth=1.1pt](4.5,-8){0.025cm}{g-1} \cnode*[linewidth=1.1pt](6,-7){0.025cm}{g-2} \cnode*[linewidth=1.1pt](7,-6){0.025cm}{g-3} \cnode*[linewidth=1.1pt](22,-6){0.025cm}{g-4}

\ncline[linewidth=1.1pt]{-}{g1}{g5} \ncline[linewidth=1.1pt]{-}{g2}{g6}   \ncline[linewidth=1.1pt]{-}{g3}{g7}  \ncline[linewidth=1.1pt]{-}{g4}{g8}
\ncline[linewidth=1.1pt]{-}{g9}{g13} \ncline[linewidth=1.1pt]{-}{g10}{g14}   \ncline[linewidth=1.1pt]{-}{g11}{g15}  \ncline[linewidth=1.1pt]{-}{g12}{g16}
\ncline[linewidth=1.1pt]{-}{g17}{a3k-2} \ncline[linewidth=1.1pt]{-}{g18}{a2k-2}   \ncline[linewidth=1.1pt]{-}{g19}{ask-2}

\ncarc[arcangle=30,linewidth=1.3pt,linecolor=black]{-}{g-3}{g3}  \ncarc[arcangle=30,linewidth=1.3pt,linecolor=black]{-}{g-2}{g2} \ncarc[arcangle=30,linewidth=1.3pt,linecolor=black]{-}{g-1}{g1}   \ncarc[arcangle=-30,linewidth=1.3pt,linecolor=black]{-}{g-4}{g4}

   \ncarc[arcangle=-10,linewidth=1.3pt,linecolor=black]{-}{g-1}{k5}
  \ncarc[arcangle=-10,linewidth=1.3pt,linecolor=black]{-}{g-2}{k3}

  \ncarc[arcangle=-30,linewidth=1.3pt,linecolor=black]{-}{a-2}{g15}

 \ncarc[arcangle=-30,linewidth=1.3pt,linecolor=black]{-}{g18}{g14}  \ncarc[arcangle=-30,linewidth=1.3pt,linecolor=black]{-}{g17}{g13}
 \ncarc[arcangle=30,linewidth=1.3pt,linecolor=black]{-}{g19}{g16}

\cnode*[linewidth=1.1pt](3,10){0.05cm}{d0} \cnode*[linewidth=1.1pt](3,11){0.05cm}{d1} \cnode*[linewidth=1.1pt](3,12){0.05cm}{d2}
\cnode*[linewidth=1.1pt](2,10){0.05cm}{d0} \cnode*[linewidth=1.1pt](2,11){0.05cm}{d1} \cnode*[linewidth=1.1pt](2,12){0.05cm}{d2}
\cnode*[linewidth=1.1pt](4,10){0.05cm}{d0} \cnode*[linewidth=1.1pt](4,11){0.05cm}{d1} \cnode*[linewidth=1.1pt](4,12){0.05cm}{d2}
\cnode*[linewidth=1.1pt](25,10){0.05cm}{d0} \cnode*[linewidth=1.1pt](25,11){0.05cm}{d1} \cnode*[linewidth=1.1pt](25,12){0.05cm}{d2}

\ncline[linewidth=1.1pt]{-}{g17}{a3k-2}    %\ncline[linewidth=1.1pt]{-}{g-2}{k3}
\ncline[linewidth=1.1pt]{-}{g-3}{k1} \ncline[linewidth=1.1pt]{-}{g-4}{ks}

\rput(10,-5.3){\small $k$}  \rput(12.8,-7){\small $3k$}  \rput(14.7,-8){\small $5k$}    \rput(20.4,-6.5){\small $(s-2)k$}
\rput(12,29.5){\small $5k-1$}    \rput(20,28.5){\small $(s-2)k-1$}   \rput(10,28.5){\small $3k-1$}

\end{pspicture}
\end{center}

\caption{A drawing of the graph $H_s(k), \ k\ge 5\/$ is odd and $s \ge 3\/$ is odd.}
\end{figure}

%\vspace{5mm} \begin{center}    {\bf Acknowledgements} \end{center} Research of the first author was supported by the  FRGS Grant (FP036-2013B).

�\end{document}